\documentclass[11pt]{elsarticle}
\usepackage{geometry}
\usepackage{graphicx}
\usepackage{verbatim}
\usepackage{array}
\usepackage{epstopdf}
\usepackage{amsmath,amsfonts,amssymb}
\usepackage{amsthm}
\usepackage[numbers]{natbib}
\usepackage{multirow}
\usepackage{setspace}
\usepackage{tabu}
\usepackage{caption}
\usepackage{multicol}
\usepackage[font=scriptsize]{subfig}
\newtheoremstyle{teorema}{\topsep}{\topsep}{}{}{}{ }{\newline}
{\fcolorbox{white}{white}{\bf \color{black}
\textsl{\thmname{#1} \thmnumber{#2}} \thmnote{\{#3\}}}}

\theoremstyle{teorema}

\begin{document}
\onehalfspacing

%%%%%%%%%%%%%%%%%%%%
%%%%%%% Front page %%%%%%%
%%%%%%%%%%%%%%%%%%%%
\begin{frontmatter}
\title{A general space-time model for combinatorial optimization problems \\ (and not only)}

\author[por]{Maria Barbati\corref{cor1}}\ead{maria.barbati@port.ac.uk}
\author[uc]{Salvatore Corrente}
\ead{salvatore.corrente@unict.it}
\author[por,uc]{Salvatore Greco}
\ead{salgreco@unict.it}

%\ead{salgreco@unict.it}
\cortext[cor1]{Corresponding author}
\address[por]{University of Portsmouth, Faculty of Business and Law, Centre of Operations Research and Logistics (CORL), Portsmouth, United Kingdom}
\address[uc]{Department of Economics and Business, University of Catania,   Catania, Italy}

%\date{}
%\maketitle

%\vspace{-1cm}

%%%%%%%%%%%%%%%%%%%%(
\begin{abstract}
We consider the problem of defining a strategy consisting of a set of facilities taking into account also the location where they have to be assigned and the time in which they have to be activated. The facilities are evaluated with respect to a set of criteria. The plan has to be devised respecting some constraints related to different aspects of the problem such as precedence restrictions due to the nature of the facilities.  Among the constraints, there are some related to the available budget. We consider also the uncertainty related to the performances of the facilities with respect to considered criteria, and plurality of stakeholders participating to the decision. The considered problem can be seen as the combination of some prototypical operations research  problems: knapsack problem, location problem and project scheduling. Indeed, the basic brick of our model is a variable $x_{ilt}$ which takes value 1 if facility $i$ is activated in location $l$ at time $t$, and $0$ otherwise. Due to the conjoint consideration of a location and a time in the decision variables, what we propose can be seen as a general space-time model for operations research problems.  We discuss how such a model permits to handle complex problems using several methodologies including multiple attribute value theory and multiobjective optimization. With respect to the latter point, without any loss of the generality, we consider the compromise programming and an interactive methodology based on the Dominance-based Rough Set Approach.  We illustrate the application of our model with a simple didactic example. \\
\end{abstract}

%evaluated with respect to a multiplicity of criteria.

\end{frontmatter}

\pagenumbering{arabic}
\journal{Omega}

%%%%%%%%%%%%%%%%%%%%%%%%%%
\section{Introduction}%%%%
%%%%%%%%%%%%%%%%%%%%%%%%%%

Operational Research (OR) has been developing around a certain number of prototypical problems such as facility location, knapsack and scheduling (for a survey see \cite{Owen}, \cite{Martello}, \cite{Hartmann}, respectively). The classical OR approach formulates these problems in terms of optimization of a well defined objective function representing the preferences of a single Decision Maker (DM) in a deterministic context.

Despite a vast number of successful applications of OR techniques, we have to admit that real world decision problems require a broader methodology than the classical OR approaches. In this perspective one can observe that in OR it is more and more common to consider a plurality of objective functions (see e.g., \cite{Deb}) taking into account preferences of a multiplicity of stakeholders (see e.g., \cite{De Gooyert}) in an uncertain environment (see e.g., \cite{Gabrel}).

We have to observe also that many real life problems present elements of more than one prototypical OR problem. Consider, for example,  the design of an urban development project in which several facilities have to be activated in different feasible locations in parallel or in a temporal sequence under some budget constraints. You can see that such a problem presents:
\begin{itemize}
	\item elements of the knapsack problem related to the facility to be selected,
	\item elements of the facility location related to the position where facilities have to be placed,
	\item elements of the scheduling problem related to the period in which the selected facilities have to be activated.
\end{itemize}

\noindent In simple words, one can say that prototypical OR problems consider only one of the following questions
\begin{itemize}
	\item ``what?'', which is the case of knapsack problems answering the question ``what items should be select?'',
	\item ``where?'', which is the case of facility location problems answering the question ``where facilities should be located?'',
	\item ``when?'', which is the case of scheduling problems answering the question ``when activities should take place?''.
\end{itemize}

\noindent Instead, complex real world decision problems consider simultaneously all the three above questions.

In the literature we can find models and methods to answer each of the single questions or pairs of them but not the combination of all three. The ``what?'' question, generated a strand of research related to the knapsack problem where one has to select which items should be inserted in a knapsack in order to optimize an objective function representing the overall profit of the items entering the knapsack while keeping the total weight of the selected items within the limited knapsack capacity. To solve the problem, in its multiobjective formulation, a good approximation of the set of solutions covering all possible trade-offs between the different objectives is identified (e.g., \cite{Captivo}, \cite{Figueira}, and \cite{MavroFig}). This strand also includes the portfolio decision problems ~\cite{Salo} in which, given a set of feasible projects evaluated on a set of criteria, the choice of the projects to insert in the portfolio is guided by the maximization of a value function of  the portfolio. Several studies have been proposed as for example in \cite{Badri}, \cite{Liesio2007}, \cite{Liesio2008} and \cite{Morton}. Often the aim is to define a new methodology to tackle these problems as in \cite{Argyris}, \cite{barbati} or \cite{Lourenco2017}.

If we are attempting to answer the ``where?'' question, we are formulating a facility location problem that consists in positioning a set of facilities in a given space. Usually the facilities have to satisfy some demand from the customer and the position of the facilities is determined on the basis of some objective functions representing the satisfaction of the demand \cite{eiseltLaporte:1995}. This strand of research is wide  and many models concerning different aspects of the problems have been proposed (for some interesting reviews  see \cite{DreznerRev}, \cite{LocationScience} or  \cite{Owen}). The classical objective is the minimisation of the sum of the distances between the users and the facilities \cite{Hakimi} but several other modifications have been proposed taking into account a huge disparity of objectives and several different constraints modeling several aspects of the problem (for a list see \cite{Farahani}).

Alternatively, considering the ``when?'' question, we are dealing with a scheduling problem consisting in defining the time in which to start different activities. The classical objective is to find a feasible schedule so that the project duration has to be minimized \cite{Habibi}. Also in this case many models have been proposed, considering different aspects of the problem such as capacity constraints \cite{Koulinas} or robustness of the problem \cite{Abbasi}.

In the literature also combinations of two of the three questions have been considered. For example, combining the ``what?'' question and the ``where?'' question, \cite{Ishizaka} selected the position for casinos in London, while \cite{Ozcan} dealt with a warehouse location selection and \cite{Tzeng} with a restaurant location selection problem. Also, \cite{Cheng} proposed a binary integer linear programming model to determine the locations for fixed investment as construction projects. Moreover, \cite{Montibeller} dealt with this type of problem considering multiple criteria and multiple stakeholders. Furthermore, the combination of the questions ``when?'' and ``where?'' has originated a flourishing strand of research related to the dynamic facility location problems (for a survey see \cite{Arabani}). In these problems the facilities are located  in each time period of a finite planning horizon \cite{Kelly}. Many modifications have been proposed as the possibility of relocating the facilities \cite{Melo2006} or the possibility of both relocating and/or changing the activated facilities and their capacities through the time in order to satisfy the customers demand \cite{Correia}. Finally, the combination of ``what?'' and ``when?'' questions is modeled in several papers in which the portfolio decision problem concerns also the timing in which each project should be developed ~\cite{Ghasemzadeh}. Again, several algorithms and methods have been proposed. For example \cite{Dickinson} and \cite{Zuluaga} introduced the interdependency of the projects, while \cite{Doerner} considered the benefits derived by the projects divided in categories. Furthermore, \cite{Ghorbani} supposed that the projects can start in some periods and continuing or not over following periods in order to maximise the benefits derived from the portfolio and the balance of the resources allocated in the different periods. Recently, \cite{Perez} modeled synergies and incompatibilities among projects and uncertainty in the parameters of the problem.

The above literature review shows that, despite very often the nature of real life problems is in between of several prototypical OR problems, there is not a general model permitting a systematic analysis of such complex problems. In view of this, we propose a general methodology permitting to handle problems that have elements of the knapsack problem, of the facility location problem and the scheduling problem at the same time. We adopt a multiobjective optimization approach to take into account a plurality of criteria as it seems natural in this type of problems. Moreover, we consider also the possibility to take into consideration uncertainty related to different potential scenarios and the presence of a plurality of stakeholders, as this can be useful in several real life contexts.

Since the problem we are handling requires answers to the basic question ``what?'' of the knapsack problem, considering also the questions ``where?'' and ``when?'', the methodology we are proposing defines a space-time model in which the activation of each facility is characterized not only in  terms of spatial coordinates typical of location problems, but also in terms of a time frame considered in scheduling problems. From the formal point of view, the basic idea is to consider variables of the type $x_{ilt}$ taking value 1 if facility $i$ is activated in location $l$ at time $t$, and 0 otherwise. Observe that our space-time model is not restricted to the above considered problems of facility location planning \cite{Cheng}, but it can be applied in other relevant situations such as, for example, a project portfolio selection \cite{Montibeller} in which, beyond the set of selected projects, it is considered the timing with which the projects have to be realized. Observe also, that, even if the three prototypical OR problems we considered are of combinatorial optimization nature, one can always relax the binary constraints permitting the decision variables to take a value on the non-negative reals. In this way our space-time model can be applied in problems that do not require combinatorial optimization, such as the typical problems of environmental planning \cite{Huang}.

The paper is organized in the following way. In Section \ref{proposed model} we present the  formulation of our space-time model. In Section \ref{example} we illustrate a didactic example for our model. In Section \ref{multiobjective} we apply two different multiobjective optimization methods to our model. In Section \ref{IncertezzaPlur} we explain how the model can be used in presence of uncertainty and plurality of stakeholders, while Section \ref{concl} concludes the paper.

\section{The proposed model}\label{proposed model}

The considered problem concerns a set of facilities $I=\{1, \ldots, i, \ldots, n\}$ to be placed in a set of feasible locations $L=\{1, \ldots, l, \ldots, m\}$ in different periods $T=\{0, 1, \ldots, t, \ldots, p\}$. Each facility is evaluated with respect to a set of criteria $J=\{1, \ldots, j, \ldots,q\}$. The evaluation of facility $i\in I$ activated in location $l \in L$ with respect to the criterion $j\in J$ is denoted by $y_{ijl} \in \mathbb{R}^+$. For the sake of simplicity, without the loss of generality, we suppose that all criteria $j\in J$ are of the gain type,  that is, the greater $y_{ijl}$, the better the evaluation of facility $i\in I$ on criterion $j \in J$ in location $l \in L$. For each period $t \in T$ a discount factor $v(t)$, with $ 0 \le v(t) \le 1$ and $v$ being a nonincreasing function of $t$, is defined in order to discount the evaluation of performances $y_{ijl}, i \in I, j \in J, l \in L$ in future periods. The values $v(t)$, $t \in T$, have to represent the intertemporal preferences of the DM. There is a vast literature on discounting and time preference (see  \cite{frederick2002time} for a survey) and, of course, among the many models proposed, that one that can be considered the more convenient with respect to the application at hand can be applied in our framework. In general, for the sake of simplicity, in the rest of the paper, when we refer to a specific discount factor $v(t)$ we consider the  model presented in \cite{samuelson1937note} and characterized by a constant interest rate $\rho$, such that $v(t)=(1+\rho)^{-t}$.  Once defined the discounted factors $v(t), t \in T$, $V_{ijlt}=y_{ijl}\cdot v(t)$ gives the value in period 0 of the performance in period $t$ of facility $i \in I$ activated in location $l \in L$ with respect to criterion $j \in J$. Here we are supposing that the benefit of facilities with respect to considered criteria does not depend on the time passed from their activation. Of course, this assumption is rather strong and can be relaxed considering a benefit depending also on the time passed from the activation, so that we have to consider an evaluation  $y_{ijlr}$ of facility $i$ activated in location $l $ with respect to the criterion $j$ after $r, r=1,\ldots,p-1$, periods from its activation. In this case, if the facility is activated in period $\tau$, the discounted value in 0 of the evaluation  $y_{ijlr}$ is given by $V_{ijl\tau r}=y_{ijlr}\cdot v(\tau+r)$. Since we need to aggregate performances on different criteria, we have to consider a weight $w_j \ge 0$ such that $w_1+\ldots+w_q=1$, for each criterion $j \in J$, in order to make homogeneous their performances and permitting their sum. Each facility $i \in I$ has also a cost $c_i\in\mathbb{R}^+$.   The available budget for each period $t \in T$ is denoted by $B_t$.

The following decision variables can be considered to define the adopted strategy:

$$x_{ilt}=
\left\{
\begin{array}{rl}
1, & \hbox{ if facility $i \in I$ is activated in location $l $ in period $t \in T-\{p\}$}; \\
0, & \hbox{ otherwise.}
\end{array} \right.
$$

For example, having a set of facilities $I=\{1,2\}$, a set of locations $L=\{1,2\}$ and a set of periods $T=\{0,1,2\}$ we have to consider the following vector decision variables:
$$\mathbf{x}=[x_{110},x_{111},x_{120},x_{121},x_{210},x_{211},x_{220},x_{221}].$$
If we have
$$x_{110}=x_{111}=x_{120}=0, x_{121}=x_{210}=1,  x_{211}=x_{220}=x_{221}=0,$$
then the adopted strategy consists in placing facility 1 in location 2 in period 1 and facility 2 in location 1 in period 0.
Observe that not all 0-1 vectors $\mathbf{x}=[x_{ilt}]$ are feasible. Indeed, some constraints have to be satisfied such 
\begin{itemize}
	\item budget constraints for which in each period $t \in T$ the expenses cannot be greater than the available budget $B_t$
\begin{equation}
\label{ConstraintBudget}\sum_{i\in I} c_i \sum_{l\in L} x_{ilt} \leq B_t, \;\; \forall  t\in T,
\end{equation}
 \item activation constraints for which each facility can be activated at most once
\begin{equation}
\label{ConstraintOnlyOneFac} \sum_{l\in L, t \in T} x_{ilt} \leq 1, \;\; \forall  i\in I.
\end{equation}
\end{itemize}
Of course, other constraints can be considered such as precedence constraints for which some facilities cannot be activated before other related facilities have already been activated. Moreover, also the budget constraints and the activation constraints can be weakened or strengthened. For example, with respect to the budget constraints, one can imagine that it is possible to lend some capital or to use the monetary return of some facility already activated. Also activation constraints can have different formulations such as no more than a fixed number of facilities of a given type can be activated.

Given a strategy $\mathbf{x}$, the benefit of criterion $j \in J$ in period $t\in T-\{0\}$ from facility $i \in I$ is obtained if $i$ has been activated not later than period $t-1$, otherwise it is null. Therefore the performance of facility $i \in I$ with respect to criterion $j \in J$ in location $l \in L$ at time $t\in T-\{0\}$ is
$$y^{IJLT}_{ijlt}(\mathbf{x})=\sum_{\tau=0}^{t-1} x_{il\tau}y_{ijl}.$$
Discounting the performance $y^{IJLT}_{ijlt}(\mathbf{x})$ we get
$$\widehat{y}^{IJLT}_{ijlt}(\mathbf{x})= y^{IJLT}_{ijlt}(\mathbf{x})v(t)=\sum_{\tau=0}^{t-1} x_{il\tau}y_{ijl}v(t).$$

Given a strategy $\mathbf{x}$, from the values $y^{IJLT}_{ijlt}(\mathbf{x})$ the following other interesting values can be obtained
\begin{itemize}
\item the global performance of the strategy $\mathbf{x}$  with respect to criterion $j \in J$ in location $l \in L$ in period $t\in T-\{0\}$ is
$$y^{JLT}_{jlt}(\mathbf{x})=\sum_{i \in I}y^{IJLT}_{ijlt}(\mathbf{x}) =\sum_{i \in I} \sum_{\tau=0}^{t-1} x_{il\tau}y_{ijl},$$
\item the overall performance of facility $i \in I$ in location $l \in L$ in period $t\in T-\{0\}$ taking into account all criteria is
$$y^{ILT}_{ilt}(\mathbf{x})=\sum_{j \in J}\sum_{\tau=0}^{t-1} w_j x_{il\tau}y_{ijl},$$
\item the performance of facility $i \in I$ with respect to criterion $j \in J$ in period $t\in T-\{0\}$ taking into account all locations is
$$y^{IJT}_{ijt}(\mathbf{x})=\sum_{l \in L} \sum_{\tau=0}^{t-1} x_{il\tau}y_{ijl},$$
\item the performance of facility $i \in I$ with respect to criterion $j \in J$ in location $l \in L$ taking into account all periods $t\in T-\{0\}$ is
$$y^{IJL}_{ijl}(\mathbf{x})=\sum_{t \in T-\{0\}} \sum_{\tau=0}^{t-1} x_{il\tau}y_{ijl},$$
\item the global performance of the strategy $\mathbf{x}$  in location $l \in L$ at time $t\in T-\{0\}$ is
$$y^{LT}_{lt}(\mathbf{x})=\sum_{i \in I} \sum_{j \in J} \sum_{\tau=0}^{t-1} w_j x_{il\tau}y_{ijl}$$
\item the performance of strategy $\mathbf{x}$ with respect to criterion $j \in J$ at time $t\in T-\{0\}$ considering all locations is
$$y^{JT}_{jt}(\mathbf{x})=\sum_{i \in I} \sum_{l \in L} \sum_{\tau=0}^{t-1} x_{il\tau}y_{ijl},$$
\item the performance of strategy $\mathbf{x}$ with respect to criterion $j \in J$ in location $l \in L$ taking into account all periods  $t\in T-\{0\}$ is
$$y^{JL}_{jl}(\mathbf{x})=\sum_{i \in I}\sum_{t \in T-\{0\}}\sum_{\tau=0}^{t-1} x_{il\tau}y_{ijl},$$
\item the performance of facility $i \in I$  in period $t\in T-\{0\}$ considering all criteria $j \in J$ and all locations $l \in L$ is
$$y^{IT}_{it}(\mathbf{x})=\sum_{j \in J}\sum_{l \in L}\sum_{\tau=0}^{t-1} w_j x_{il\tau}y_{ijl},$$
\item the performance of facility $i \in I$  in location $l \in L$ considering all criteria $j \in J$ and all periods $t\in T-\{0\}$ is
$$y^{IL}_{il}(\mathbf{x})=\sum_{j \in J}\sum_{t \in T-\{0\}}\sum_{\tau=0}^{t-1} w_j x_{il\tau}y_{ijl},$$
\item the performance of facility $i \in I$ with respect to criterion $j \in J$ considering all locations  $l \in L$ and all periods from $t\in T-\{0\}$ is
$$y^{IJ}_{ij}(\mathbf{x})=\sum_{l \in L} \sum_{t \in T-\{0\}} \sum_{\tau=0}^{t-1} x_{il\tau}y_{ijl},$$
\item the overall performance in period $t\in T-\{0\}$ considering all facilities $i \in I$, all criteria $j \in J$ and all  locations $l \in L$ is
$$y^{T}_{t}(\mathbf{x})=\sum_{i \in I}\sum_{j \in J}\sum_{l \in L}\sum_{\tau=0}^{t-1} w_j x_{il\tau}y_{ijl},$$
\item the overall performances of strategy $\mathbf{x}$  in location $l \in L$ considering all criteria $j \in J$ and all periods $t\in T-\{0\}$ is
$$y^{L}_{l}(\mathbf{x})=\sum_{i \in I}\sum_{j \in J}\sum_{t \in T-\{0\}}\sum_{\tau=0}^{t-1} w_j x_{il\tau}y_{ijl},$$
\item the overall performance with respect to criterion $j \in J$ considering all facilities $i \in I$, all locations $l \in L$ and all periods $t\in T-\{0\}$ is
$$y^{J}_{j}(\mathbf{x})=\sum_{i \in I}\sum_{l \in L}\sum_{t \in T-\{0\}}\sum_{\tau=0}^{t-1} x_{il\tau}y_{ijl},$$
\item the performance of facility $i \in I$ considering all criteria  $j \in J$, all locations $l \in L$ and all periods $t\in T-\{0\}$ is
$$y^{I}_{i}(\mathbf{x})=\sum_{j \in J}\sum_{l \in L}\sum_{t \in T-\{0\}}\sum_{\tau=0}^{t-1} w_j x_{il\tau}y_{ijl},$$
\item the overall performance of strategy $\mathbf{x}$ taking into account all facilities $i \in I$, all criteria $j \in J$, all locations $l \in L$ and all periods $t\in T-\{0\}$ is
$$y(\mathbf{x})=\sum_{i \in I}\sum_{j \in J}\sum_{l \in L}\sum_{t \in T-\{0\}}\sum_{\tau=0}^{t-1} w_j x_{il\tau}y_{ijl}.$$
\end{itemize}

Let us point out that all the above performances can be discounted. For example the discounted value at time 0 of $y^{JLT}_{jlt}(\mathbf{x})$ is given by
$$\widehat{y}^{JLT}_{jlt}(\mathbf{x})=y^{JLT}_{jlt}(\mathbf{x})v(t)=\sum_{i \in I} \sum_{\tau=0}^{t-1} x_{il\tau}y_{ijl}v(t).$$
We shall denote with $\widehat{y}^{sets}_{indices}(\mathbf{x})$ the discounted value of the corresponding non-discounted performance $y^{sets}_{indices}(\mathbf{x})$, so that $\widehat{y}^{JLT}_{jlt}(\mathbf{x})$ is the discounted value of $y^{JLT}_{jlt}(\mathbf{x})$, $\widehat{y}^{ILT}_{ilt}(\mathbf{x})$ is the discounted value of $y^{ILT}_{ilt}(\mathbf{x})$, and so on.

In the first instance, the problem is to define the strategy $\mathbf{x}$ giving the maximum overall discounted performance $\widehat{y}(\mathbf{x})$ subject to the constraints of the problem such as the budget constraints and the activation constraints.

However, the above model permits to take into account a great plurality of performances $y^{sets}_{indices}(\mathbf{x})$ and $\widehat{y}^{sets}_{indices}(\mathbf{x})$  constituting a rich dashboard that can be very meaningful for the DM. In fact, the DM can fix some constraints in terms of minimal requirements of performances $y^{sets}_{indices}(\mathbf{x})$ and $\widehat{y}^{sets}_{indices}(\mathbf{x})$. More in general, we can handle the whole model in terms of multiobjective optimization of  performances $y^{sets}_{indices}(\mathbf{x})$ and $\widehat{y}^{sets}_{indices}(\mathbf{x})$.  We shall explore this possibility in Section \ref{multiobjective}.

\subsection{A possible extension with continuous variables}

Our model can work also when the variables $x_{ilt}$ are defined in $\mathbf{R^+}$. In this case the variables can be defined as the amount of budget that has been allocated to facility of type $i$ in location $l$ at period $t$. In this case binary constraints (\ref{ConstraintOnlyOneFac}) must not be considered, while to the original budget constraint (\ref{ConstraintBudget}), we can add additional budget constraints. In particular,  we can define
\begin{itemize}
\item  $B^\leq_{it}$ as the maximum budget to be allocated to facility $i \in I$ in  period $t \in T$,
\item $B^\geq_{it}$ as the minimum budget to be allocated to facility $i \in I$ in  period $t \in T$,
\item $B^\leq_{lt}$ as the maximum budget to be allocated to location $l \in L$ in  period $t \in T$,
\item $B^\geq_{lt}$ as the minimum budget to be allocated to to facility $l \in L$ in  period $t \in T$,
\end {itemize}

so that for each of the above quantities we can define the additional constraints:

\begin{itemize}
	\item in period $t \in T$, not more than the maximum budget  $B^\leq_{it}$ can be allocated to facility $i \in I$:
\begin{equation}
\label{ConstraintBudgetFaciliy}\sum_{l\in L}x_{ilt} \leq B^\leq_{it}, 
\end{equation}
 \item in  period $t \in T$, not less than the minimum budget $B^\geq_{it}$ must be allocated to facility $i \in I$:
\begin{equation}
\label{ConstraintBudgetLocation} \sum_{l\in L}  x_{ilt} \geq B^\geq_{it}, 
\end{equation}
\item in  period $t \in T$, not more than the maximum budget  $B^\leq_{lt}$ can be allocated to location $l \in L$:
\begin{equation}
\label{ConstraintBudgetFaciliy}\sum_{l\in L}  x_{ilt} \leq B^\leq_{lt}, 
\end{equation}
 \item in  period $t \in T$, not less than the minimum budget $B^\geq_{lt}$ must be allocated to location $l \in L$:
\begin{equation}
\label{ConstraintBudgetLocation} \sum_{l\in L}  x_{ilt} \geq B^\geq_{it}.
\end{equation}
\end{itemize}

Those budgets, and the associated constraints, are not necessarily defined for all the facilities $i \in I$ and for all the locations $l \in L$. To ensure the feasibility of the model, it should also be verified that for each $i \in I$ then $B^\leq_{it}\leq B_t$  and for each $l \in L$ then $B^\leq_{lt}\leq B_t$.
Let us underline that also for the continuous case we can handle the whole model in terms of multiobjective optimization of the performances $y^{sets}_{indices}(\mathbf{x})$ and $\widehat{y}^{sets}_{indices}(\mathbf{x})$.

%%%%%%%%%%%%%%%%%%%%%%%%%%%%%%%%%%%%%%%%%%%%%%%%%%%%%%%%%%%%%
\section{Illustrative example}\label{example}
%%%%%%%%%%%%%%%%%%%%%%%%%%%%%%%%%%%%%%%%%%%%%%%%%%%%%%%%%%%%%

We illustrate the proposed model with the following hypothetical decision problem. Let us suppose that a council is expected to decide which public interest facilities should be activated in the next 5 years, choosing between two possible locations available for each of them. In particular, we consider an example involving the following eight desirable facilities $I=\{1, \ldots, 8\}$:

\begin{itemize}
\item School, $i=1$,
\item Leisure Centre, $i=2$,
\item Council Offices, $i=3$,
\item Recycling Centre, $i=4$,
\item Start Up Incubator, $i=5$,
\item Healthcare Centre, $i=6$,
\item Community Centre, $i=7$,
\item Social Housing, $i=8$.
\end{itemize}

\noindent evaluated in terms of the following three criteria $J=\{1, \ldots, 3\}$:

\begin{itemize}
\item Economic impact, $j=1$,
\item Social impact, $j=2$,
\item Environmental impact, $j=3$,.
\end{itemize}

\begin{table}
\centering\caption{Evaluations on the three criteria in each location and associated costs for the eight facilities considered in the illustrative example.}\label{Facilities_evaluation}

\begin{tabu}{|l||c|c|c|c|c|c||c|}
\hline
   \multirow{2}{*}{Facilities}&   \multicolumn{2}{|c|}{$Economic Impact$} &  \multicolumn{2}{|c|}{$Social Impact$} &   \multicolumn{2}{|c|}{$Environmental Impact$} &\multirow{2}{*}{Cost}  \\\cline{2-7}

    	 &	North	&	South	&	North	&	South	&	North       	&	South& 	\\
      \hline
    	School &	21	&	23	&	90	&	80	&	23	&	32&200	\\
 Leisure Centre&  36&	46&	59&	72&	36&	34&	300\\
 Council Offices& 18&	20&	22&	30&	21&	26&	150\\
Recycling Centre & 60&	65&	71&	60&	90&	88&	100\\
Start Up Incubator & 80&	82&	12&	12&	15&	12&	150\\
Healthcare  Centre & 20&	18&	19&	19&	45&	59&	200\\
Community Centre&35&	31&	56&	48&	33&	40&	100\\
Social Housing &12&	21&	69&	73&	18&	17&	250\\

     \hline
\end{tabu}
\end{table}

\begin{table}
\centering\caption{Budget available in each period}\label{Budget}

\begin{tabu}{|l||c|}
\hline
Year	&	Budget		\\
\hline
Start &	400		\\
First Year&  100\\
Second Year& 200\\
Third Year & 200\\
Fourth Year & 150\\
\hline
\end{tabu}
\end{table}

We suppose to have two different locations $L=\{1, 2\}$ in which the facilities can be positioned, named North ($l=1$) and South ($l=2$).
For the sake of the simplicity, we give an evaluation of each facility on each criterion and for each location on a scale [0,100] (see Table~\ref{Facilities_evaluation}). We assume that the evaluation does not depend on the period. Note, however, that our model can deal also with evaluations that change through the time and with any type of quantitative evaluations. Moreover, each facility has an associated opening cost (in thousand Euro) which is also reported in Table~\ref{Facilities_evaluation}. The available budget (in thousand Euro) is given for each period and detailed in Table~\ref{Budget}. In addition, the interest rate is supposed to be equal to $0.1$ for all the periods.
The council is setting up the plans for the next 5 years $T=\{0, 1, \ldots, 5\}$  deciding which investments pursuit. We define a weight for each criterion, and in particular $w_1=0.5$ for the economic impact, $w_2=0.3$ for the social impact and $w_3=0.2$ for the environmental impact.

Using the commercial software CPLEX v.12.1, we find the vector $\mathbf{x}$,  that maximises the objective function $\widehat{y}(\mathbf{x})$ subject to the budget constraint. We also suppose that each facility can be activated only once, e.g., each facility cannot be activated in two different locations and in two different periods.

\begin{figure}
\centering
\captionsetup[subfigure]{labelformat=empty}
\subfloat {\includegraphics[width = 6in]{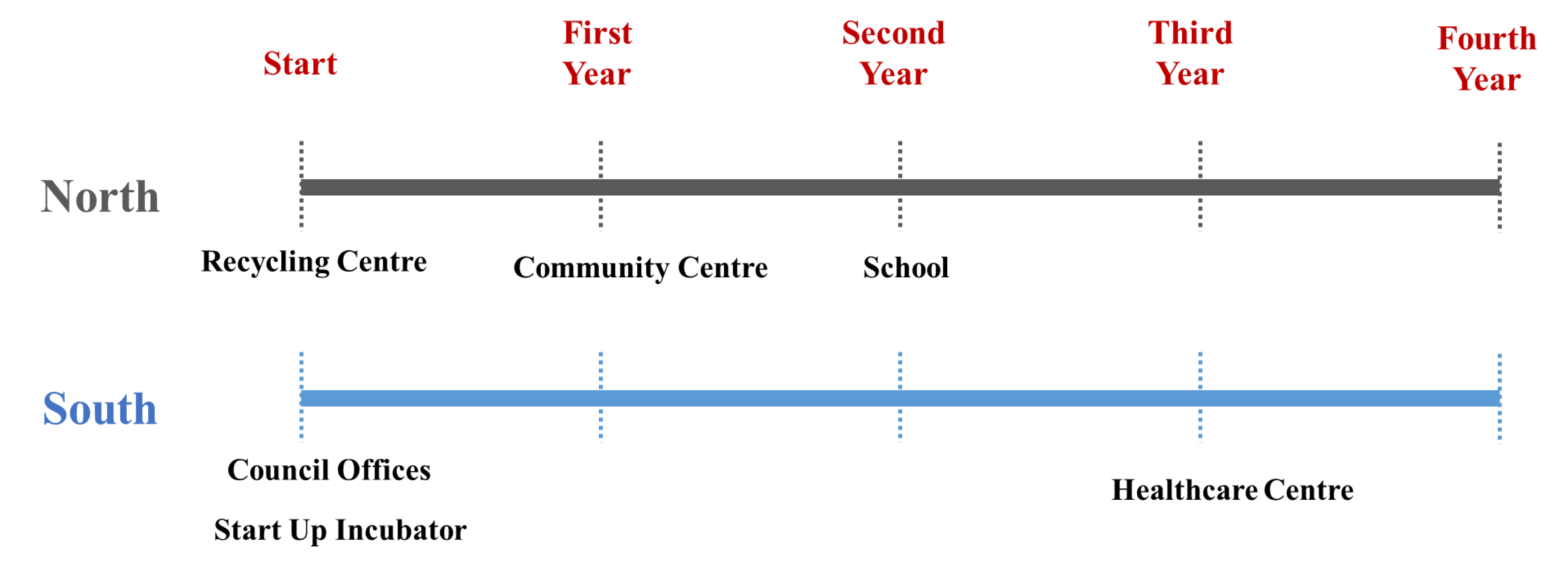}}
\caption{Optimal solution obtained by maxixization of the overall performance.}
\label{UtilitySolution}
\end{figure}

We obtain the following decision variables equals to 1: $x_{112},x_{320},x_{410},x_{520},x_{623},x_{711},$  meaning that:
\begin{itemize}
  \item The facility \emph{School} is scheduled to be activated in location \emph{North} at the beginning of the \emph{second} year;
  \item The facility \emph{Council Offices} is scheduled to be activated in location \emph{South} at the beginning of the \emph{start} year;
  \item The facility \emph{Recycling Centre} is scheduled to be activated in location \emph{North} at the beginning of the \emph{start} year;
  \item The facility \emph{Start Up Incubator} is scheduled to be activated in location \emph{South} at the beginning of the \emph{start} year;
  \item The facility \emph{Healthcare Centre} is scheduled to be activated in location \emph{South} at the beginning of the \emph{third} year;
  \item The facility \emph{Community Centre} is scheduled to be activated in location \emph{North} at the beginning of the \emph{first} year.
\end{itemize}

\noindent The other facilities (Leisure Centre and Social Housing) have not been activated given the available budget constraint. The optimal strategy is reported in Figure \ref{UtilitySolution}.

\begin{figure}
\centering
\captionsetup[subfigure]{labelformat=empty}
\subfloat {\includegraphics[width = 7in]{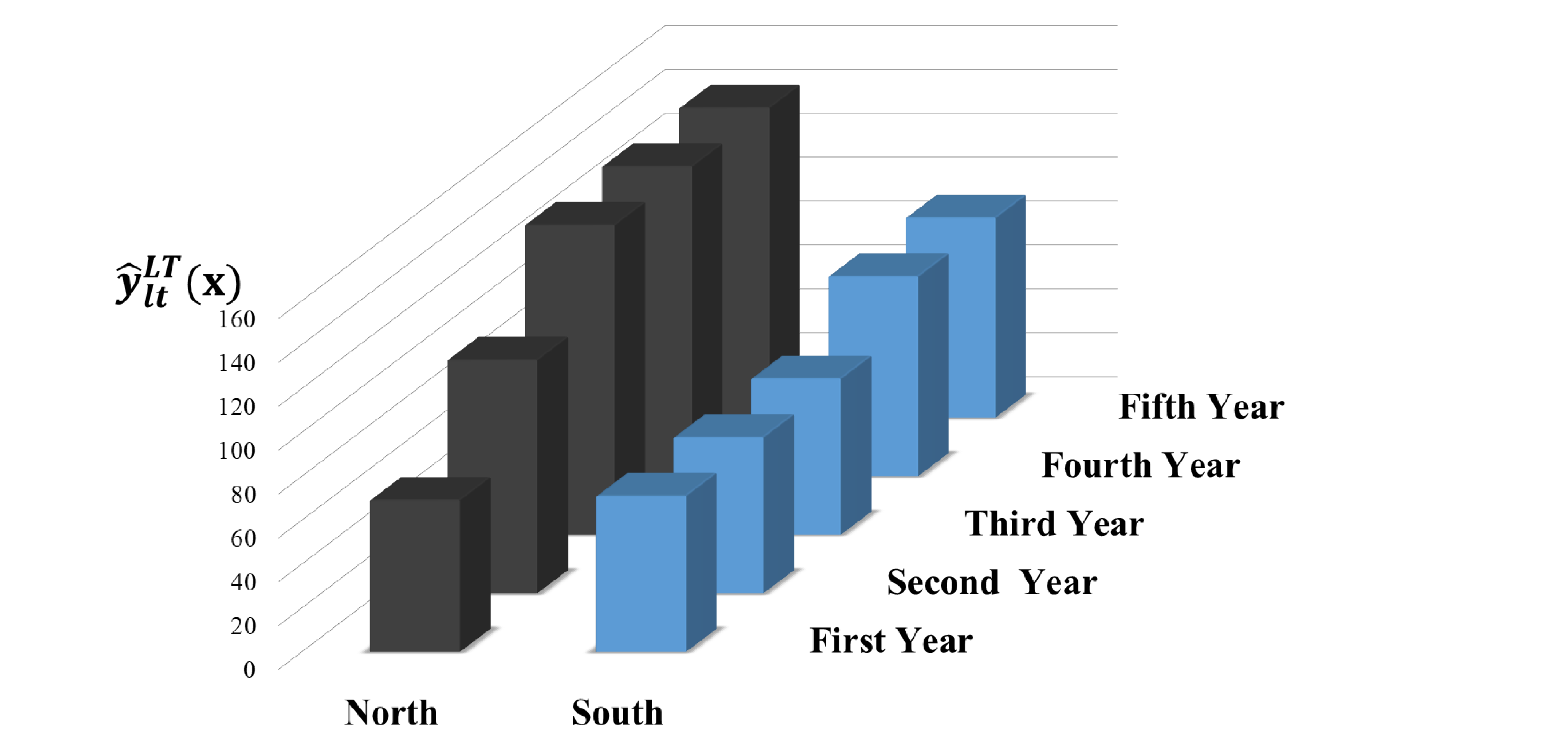}}
\caption{Distribution along the time of the performances for the optimal strategy  in the two different locations.}
\label{LocationTempo}
\end{figure}

\begin{figure}
\centering
\captionsetup[subfigure]{labelformat=empty}
\subfloat {\includegraphics[width = 7in]{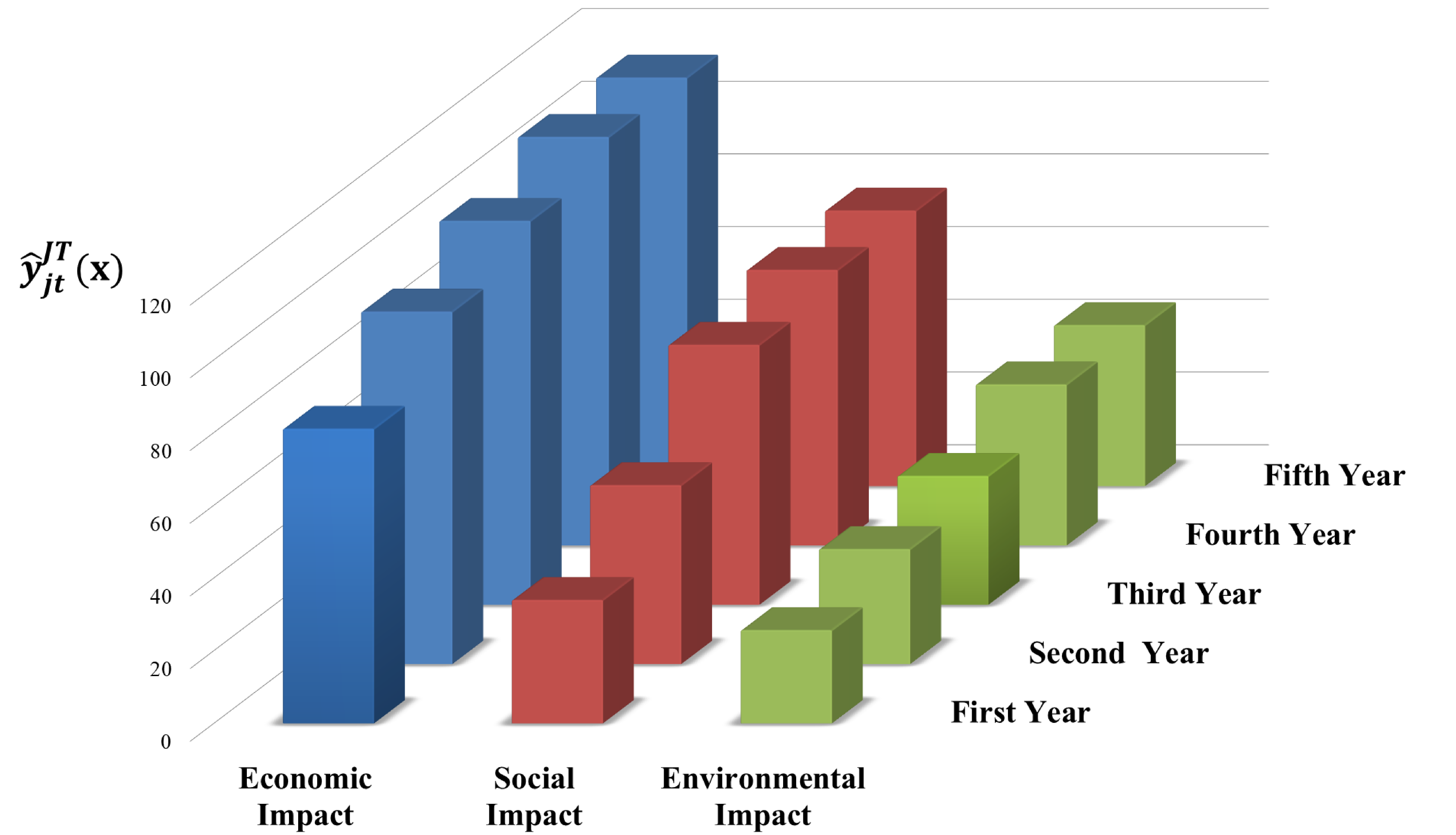}}
\caption{Time distribution of the performances for the optimal strategy with respect to each criterion.}
\label{AttributiTempo}
\end{figure}

The performances $\widehat{y}^{sets}_{indices}(\mathbf{x})$ that we have described before can be summarized in a series of graphs. These graphs can help understanding the solution and especially can help the DM to visualize the performance corresponding to the optimal solution \cite{silvamei}. Indeed, these charts can be used to compare potential Pareto solutions in a multiobjective context, supporting the intuition of the DM, and making the model more appealing even for high level managers often inhibited from adopting more sophisticated and complex decision support models~\cite{Ghasemzadeh}. For the sake of space we present the most representative charts.

First, let us show in Figure \ref{LocationTempo} the performance of the strategy suggested to the council (i.e., the optimal solution to our time - space model) in the two locations North and South. In this case we are recording the $\widehat{y}^{LT}_{lt}(\mathbf{x})$ (on the $y-axis$) in each period and in each location. It is possible to note that the performance assume a bigger value in the North than in the South. Also, there is an increasing of the performance in the time at a greater pace in location North than in location South. Note that the performances are defined excluding period $0$, that represents the start of our planning horizon. While at $t=0$ we can define decision variables, the performance of the adopted plan will be evaluated only at the beginning of the first year. Moreover, the performances are discounted so that we can compare the contribution of performances obtained in different periods.

\begin{figure}
\centering
\captionsetup[subfigure]{labelformat=empty}
\subfloat {\includegraphics[width =7.5 in]{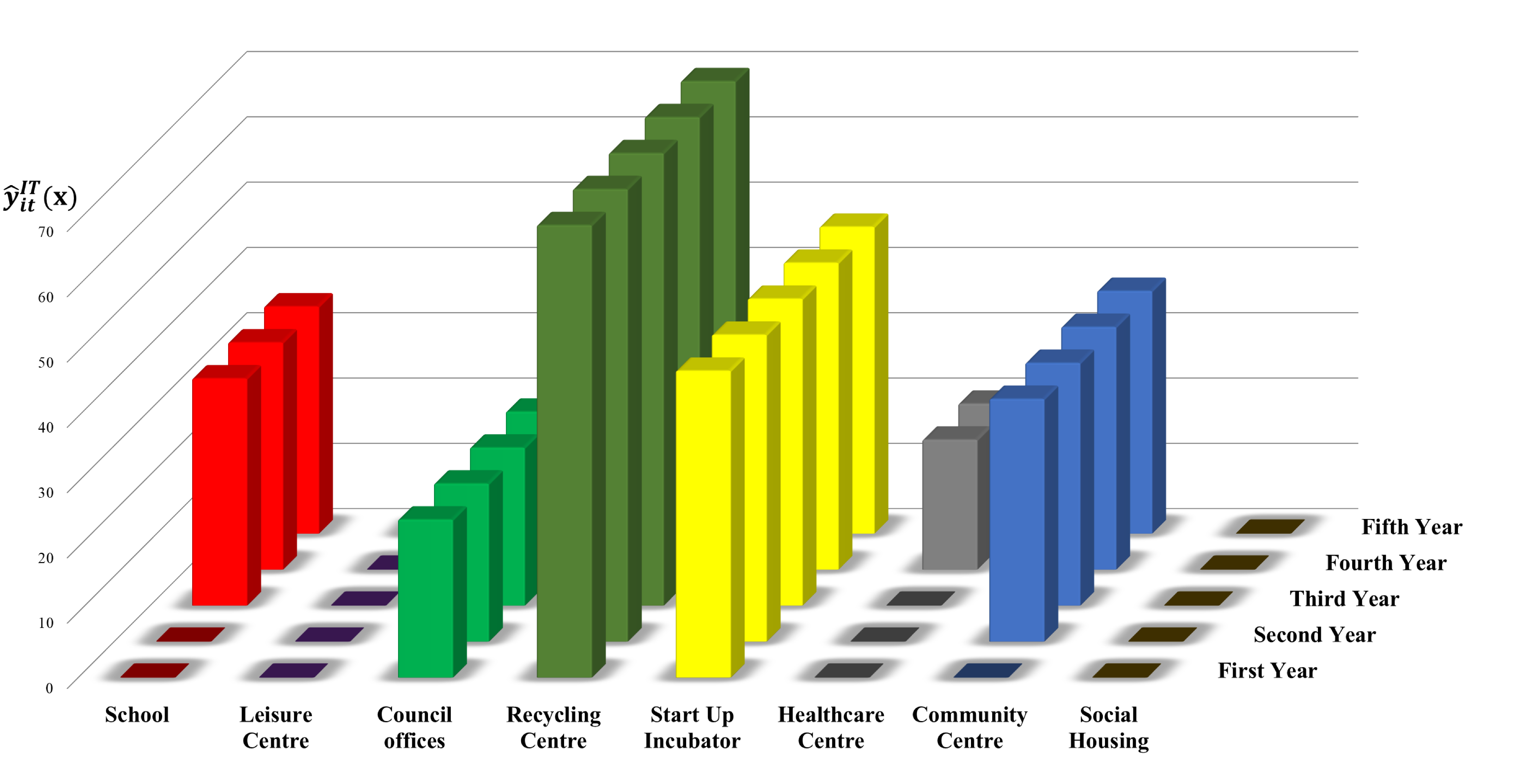}}
\caption{Time distribution of the performances for the facilities in the optimal strategy.}
\label{FacilitiesTempo}
\end{figure}

\begin{figure}
\centering
\captionsetup[subfigure]{labelformat=empty}
\subfloat {\includegraphics[width = 8in]{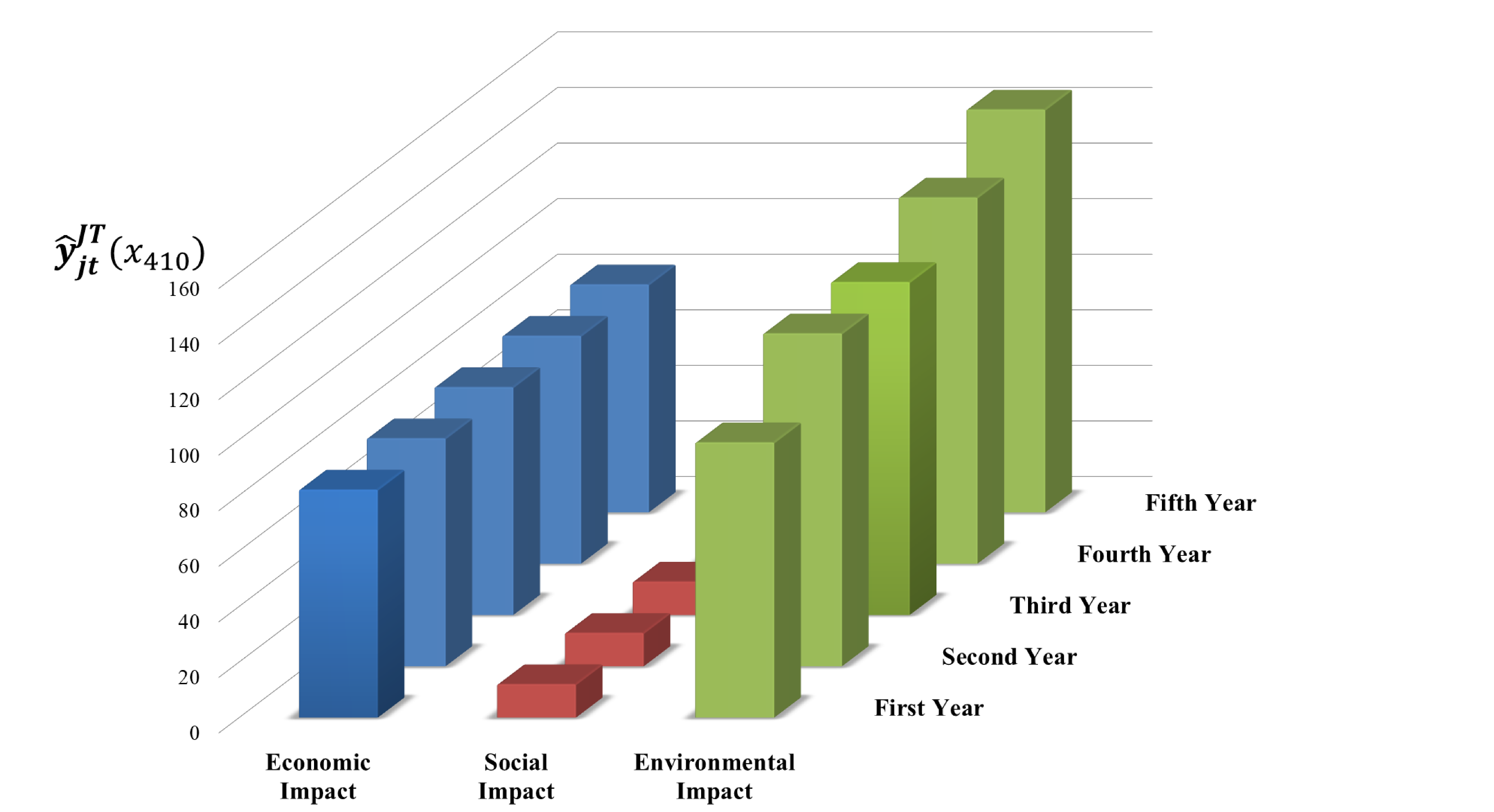}}
\caption{Time distribution of the performances for the facility \emph{Recycling Centre}  respect to each criterion.}
\label{Recycling}
\end{figure}

\begin{figure}
\centering
\captionsetup[subfigure]{labelformat=empty}
\subfloat {\includegraphics[width =6 in]{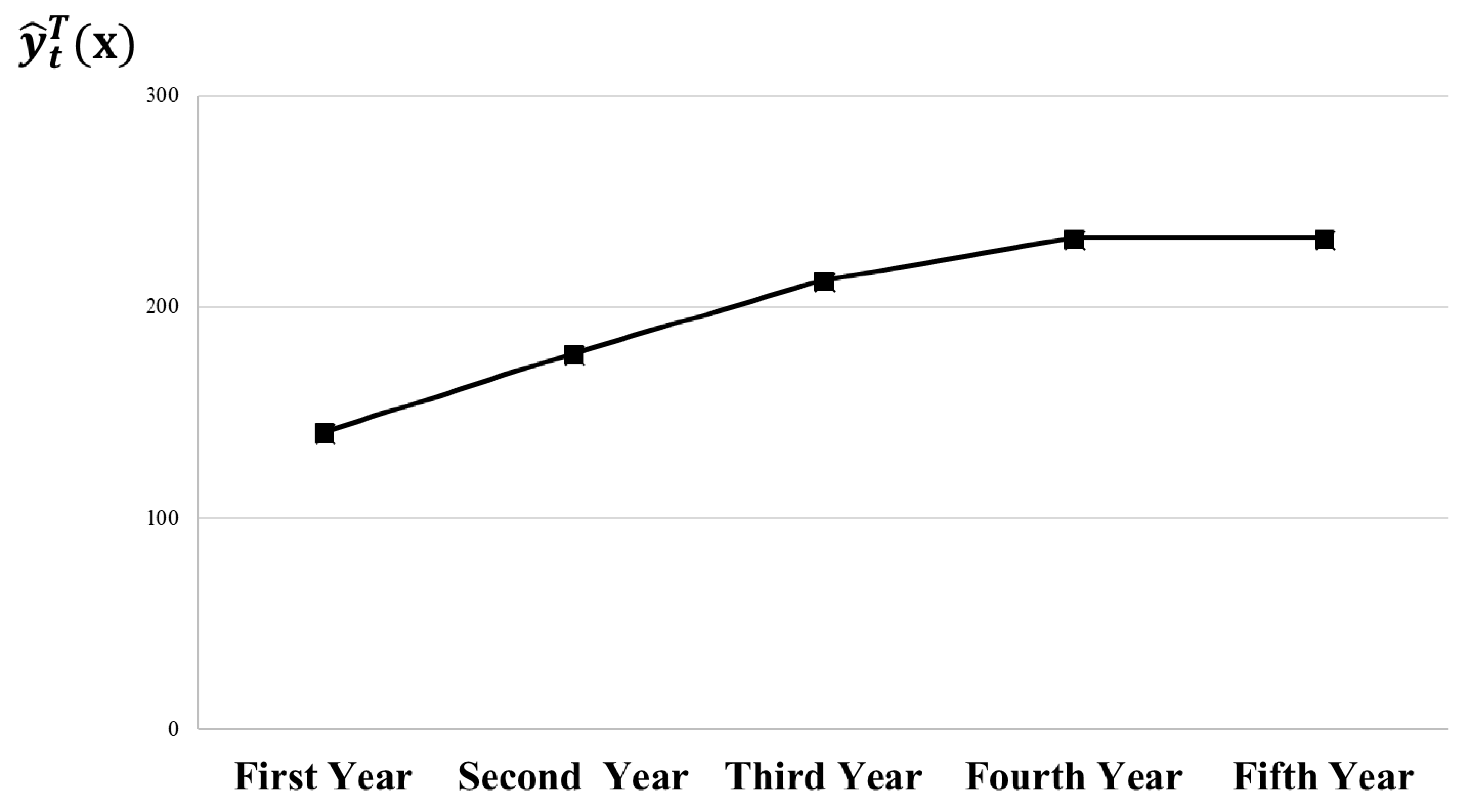}}
\caption{Time distribution of the overall performances.}
\label{Complessivo}
\end{figure}

In Figure \ref{AttributiTempo} we report the performance $\widehat{y}^{JT}_{jt}(\mathbf{x})$ of the optimal strategy with respect to each criterion and through the time.  We can see that Economic impact has a greater importance for the solution given the highest bars and its bigger increase through the time. Note that the performances in this chart have not been weighted. This allows a neat comparison without the influence of particular weights adopted.
In Figure \ref{FacilitiesTempo} we summarize the performance of each activated facility in the optimal solution  through the time. We are representing $\widehat{y}^{IT}_{it}(\mathbf{x})$, i.e., the performance of each activated facility through the time. Indeed, the  \emph{Recycling Centre} is contributing more than the other facilities. The DM could be interested in detailing the contribution of this facilities for each criterion, in the location North where it has been activated (see Figure \ref{Recycling}). The biggest contribution is provided by the criterion Environmental Impact. In this case we are reporting the $\widehat{y}^{JT}_{jt}(x_{410})$.

Finally, in Figure \ref{Complessivo} we can summarize the overall performance of the optimal strategy provided by the activated facilities through the time, indicating on the $y-axis$ the $\widehat{y}^{T}_{t}(\mathbf{x})$. This graph can help DMs to visualize the increase through the time of the contribution of all facilities, for all criteria and for each location. For this solution we can highlight that the increasing has a similar pace for the first four years while is less strong in the final year.

\subsection{Illustrative Example: Continuous case}
For the continuous case, the values of $B^\geq_{it}$, $B^\leq_{it}$, $B^\geq_{lt}$ and $B^\leq_{lt}$ are reported in Tables \ref{BudgetFac} and \ref{BudgetLoc1}, respectively. Those budgets are not defined in all the cases (when no value is defined for the budget, and so no associated constraint is defined, the symbol ``$-$" is reported in the Tables). These values are used in the formulation of constraints  \ref{ConstraintBudgetFaciliy} and \ref{ConstraintBudgetLocation}. Adopting, as before, the weighted approach, we obtain the  variables different from 0.

\begin{table}
\centering\caption{Budget available for each type of facility in each period}\label{BudgetFac}
\begin{tabu}{|l||c|c|c|c|c|c|c|c|c|c|}
\hline
   \multirow{2}{*}{Facilities}&   \multicolumn{2}{|c|}{Start} &  \multicolumn{2}{|c|}{First Year} &   \multicolumn{2}{|c|}{Second Year} &\multicolumn{2}{|c|}{Third Year} &   \multicolumn{2}{|c|}{Fourth Year} \\\cline{2-11}

    	 &	$B^\leq_{i0}$	&	$B^\geq_{i0}$	&	$B^\leq_{i1}$	&	$B^\geq_{i1}$ &	$B^\leq_{i2}$	&	$B^\geq_{i2}$ &	$B^\leq_{i3}$	&	$B^\geq_{i3}$ &	$B^\leq_{i4}$	&	$B^\geq_{i4}$ 	\\
      \hline
School 		&	-&10&	50&-&	-&5&	20&-&	-&-	\\
Leisure Centre  		&	70&20&	150&5&	10&10&	-&-&	-&16	\\
Council Offices 		&	-&-&	-&3&	8&8&	-&-&	-&-	\\
Recycling Centre  		&	32&16&	-&-&	-&-&	2&2&	260&5	\\
Start Up Incubator  		&	-&-&	70&5&	140&10&	-&2&	14&-	\\
Healthcare  Centre  		&	-&8&	-&-&	-&4&	-&-&	-&2	\\
Community Centre		&	30&-&	-&5&	-&-&	10&1&	16&14	\\
Social Housing		&	-&16&	60&-&	180&10&	-&-&	-&-	\\
     \hline
\end{tabu}
\end{table}

\begin{table}
\centering\caption{Budget available for each location in each period}\label{BudgetLoc1}
\begin{tabu}{|l||c|c|c|c|c|c|c|c|c|c|}
\hline
   \multirow{2}{*}{Locations}&   \multicolumn{2}{|c|}{Start} &  \multicolumn{2}{|c|}{First Year} &   \multicolumn{2}{|c|}{Second Year} &\multicolumn{2}{|c|}{Third Year} &   \multicolumn{2}{|c|}{Fourth Year} \\\cline{2-11}

    	 &	$B^\leq_{l0}$	&	$B^\geq_{l0}$	&	$B^\leq_{l1}$	&	$B^\geq_{l1}$ &	$B^\leq_{l2}$	&	$B^\geq_{l2}$ &	$B^\leq_{l3}$	&	$B^\geq_{l3}$ &	$B^\leq_{l4}$	&	$B^\geq_{l4}$ 	\\
      \hline
North	&	65&6&	-&2&	-&-&	4&4&	13&10	\\
South	&	21&3&	10&7&	6&6&	20&5&	-&-	\\
     \hline
\end{tabu}
\end{table}

\begin{table}
\centering\caption{Budget available for each location in each period}\label{BudgetLoc}
\begin{tabu}{|c|c|}
\hline
Decision Variables& Values \\
\hline
$x_{110}$	&	10	\\
$x_{112}$	&	5	\\
$x_{220}$	&	20	\\
$x_{221}$	&	5	\\
$x_{222}$	&	10	\\
$x_{223}$   &	195	\\
$x_{224}$	&	16	\\
$x_{321}$	&	3	\\
$x_{322}$	&	8	\\
$x_{410}$	&	32	\\
$x_{411}$	&	82	\\
$x_{412}$	&	153	\\
$x_{413}$	&	2	\\
$x_{414}$	&	118	\\
$x_{511}$	&	3	\\
$x_{513}$	&	1	\\
$x_{521}$	&	2	\\
$x_{522}$	&	10	\\
$x_{523}$	&	1	\\
$x_{610}$	&	8	\\
$x_{622}$	&	4	\\
$x_{624}$	&	2	\\
$x_{711}$	&	5	\\
$x_{713}$	&	1	\\
$x_{714}$	&	14	\\
$x_{810}$	&	15	\\
$x_{820}$	&	1	\\
$x_{822}$	&	10	\\
\hline
\end{tabu}
\end{table}

We obtain the temporal distribution of the budget between facilities and locations shown in Table \ref{BudgetLoc}. For example, the variable $x_{110}=10$, means that 10 is the budget allocated for the activation of \emph{School} in location North at the start of the planning period; the variable $x_{112} =5$ means that 5 is the budget allocated to the activation of \emph{School} in the first year in location South, and so on.

Graphs and charts analogous to those ones reported for the combinatorial model  can be provided also in this case.

%%%%%%%%%%%%%%%%%%%%%%%%%%%%%%%%%%%%%%%%%%%%%%%%%%%%%%%%%%%%%%%%
\section{Multiobjective methodologies for the space - time model}\label{multiobjective}
%%%%%%%%%%%%%%%%%%%%%%%%%%%%%%%%%%%%%%%%%%%%%%%%%%%%%%%%%%%%%%%

Several algorithms, mainly exact, have been provided in the literature to find solutions to multi-objective 0-1 linear programming problems (for a review, see \cite{EhrgottRev}).  When dealing with small problem instances, some algorithms can look for an approximation of the whole set of efficient solutions. These include the branch and bound algorithms~\cite{Przybylski} or the $\epsilon$ constraint method \cite{Cohon, Mavrotasepsilon}. Some interactive algorithms integrate optimization procedures (i.e., \cite{Alves, Argyris, Mavrotas}) with the aim of singling out the set (possibly a singleton) of the most preferred solutions for the DM. In the same perspective, other methods suggest the adoption of a linear value approach (see, e.g., \cite{Salo}) or the use of a goal programming procedure~\cite{Jones}. In what follows, we illustrate how our space - time model can be handled with two multiobjective methodologies. First we consider a classical approach called Compromise Programming (CP)  \cite{Romero} adopted to solve several multiobjectvie optimization models. The second approach is more recent and takes into account the preferences of the DM using an interactive procedure. It has been proposed by \cite{gms_IMO} and applied to portfolio decision problems in \cite{barbati}.

\subsection{Compromise Programming}

In a CP approach the aim is to minimize the maximum deviation from the ideal point, i.e., the point with the best evaluation. For our model we characterize three types of CP approaches considering three different ideal points.

First, in the Compromise Programming for Location \textbf{(CPL)}  we characterize our target as
the vector $\widehat{y}^{L^\ast}=[\widehat{y}_l^{L^\ast}]$ where, for each $l \in L$, $\widehat{y}_l^{L*}$ represents the best actualized performance that can be attained by location $l$, that is
$$\widehat{y}^{L^\ast}_l=\max_{\mathbf{x}}\widehat{y}^{L}_{l}(\mathbf{x})=\sum_{i \in I}\sum_{j \in J}\sum_{t \in T-\{0\}}\sum_{\tau=0}^{t-1} w_j x_{il\tau}y_{ijl}v(\tau)$$.

Different metrics can be adopted in order to define the closeness of the obtained strategy to the ideal point. Following \cite{Drezner2}, in order to get a balanced solution, we minimize the maximum relative deviation $\Delta_l^L(\mathbf{x})$, on the set of locations $l \in L$, defined as
$$
\Delta_l^L(\mathbf{x})=  \frac{\widehat{y}^{L^\ast}_l - \widehat{y}^{L}_{l}(\mathbf{x})}{\widehat{y}^{L^\ast}_l}.
$$
Then, the distance of the strategy  $ \mathbf{x}$ from the ideal point is $\Delta^L(\mathbf{x})=\max_{l\in L} \Delta_l^L(\mathbf{x})$.
This optimisation strategy could suit several DMs. In our example the council could be interested in attempting to minimize the differences among the locations so that the optimal solution is $\mathbf{x}^\ast=\arg\min\Delta^L(\mathbf{x})$.

Second, we  specify what we call the Compromise Programming for Objectives \textbf{(CPO)} where  the target is
the vector $\widehat{y}^{J^*}=[\hat{y}_j^{J^*}]$ where for each $j \in J$, $\widehat{y}_j^{J^*}$ represents the best actualized performance that can be attained on criterion $j$, that is
$$\widehat{y}^{J^\ast}_{j}=\max_{\mathbf{x}}\widehat{y}^{J}_{j}(\mathbf{x})=\sum_{i \in I}\sum_{l \in L}\sum_{t \in T-\{0\}}\sum_{\tau=0}^{t-1} w_j x_{il\tau}y_{ijl}v(\tau)$$. Analogously to the previous case, we shall minimize the  maximum relative deviation $\Delta_j^J(\mathbf{x})$, on the set of criteria $j \in J$, defined as
$$\Delta_j^J(\mathbf{x})=  \frac{\widehat{y}^{J^\ast}_{j} - \widehat{y}^{J}_{j}(\mathbf{x})}{\widehat{y}^{J^\ast}_{j}}.$$

Then, the distance of the strategy  $x \in \mathbf{x}$ from the ideal point is $\Delta^J(\mathbf{x})=\max_{j\in J} \Delta_j^J(\mathbf{x})$.
DMs adopting such an optimization strategy would like to balance the importance of all the criteria so that the optimal solution is $\mathbf{x}^\ast=\arg\min\Delta^J(\mathbf{x})$.

Lastly, we can define what we call Compromise Programming for Objectives and Location \textbf{(CPOL)} where  our target is $$\widehat{y}^{JL^\ast}_{jl}=\max_{\mathbf{x}}\widehat{y}^{JL}_{jl}(\mathbf{x})=\sum_{i \in I}\sum_{t \in T-\{0\}}\sum_{\tau=0}^{t-1} w_j x_{il\tau}y_{ijl}v(\tau)$$. Again, we shall minimize the maximum relative deviation $\Delta_{jl}^{JL}(\mathbf{x})$, on the set of criteria $j \in J$ and on the set of locations $l \in L$, defined as

$$\Delta_{jl}^{JL}(\mathbf{x})=  \frac{\widehat{y}^{JL^\ast}_{jl} - \widehat{y}^{JL}_{jl}(\mathbf{x})}{\widehat{y}^{JL^\ast}_{jl}}.$$

Then, the distance of the strategy  $\mathbf{x}$ from the ideal point is $\Delta^{JL}(\mathbf{x})=\max_{j\in J, l \in L} \Delta_{jl}^{JL}(\mathbf{x})$.
This last case is a combination of the first two compromise optimization approaches and attempts to balance the differences from the ideal points for both the criteria and the locations $\mathbf{x}^\ast=\arg\min\Delta^{JL}(\mathbf{x})$.

\vspace{5mm}
\noindent \emph{Illustrative Example: Compromise Programming}

We apply the three compromise optimization approaches described above to our illustrative example introduced in Section \ref{example}. We obtain the optimal compromised strategies reported in Table \ref{solcompr}. For each CP approach we noted the location and the period in which a facility has been activated; the symbol ``$-$" means that a facility has not been activated.

\begin{table}
\centering\caption{Position and activation period of the facilities in the best solutions obtained with the different CP approaches.}\label{solcompr}

\begin{tabu}{|l||c|c||c|c||c|c|}
\hline
   \multirow{2}{*}{Facilities}&   \multicolumn{2}{|c||}{$\textbf{CPL}$} &  \multicolumn{2}{|c||}{$\textbf{CPO}$} &   \multicolumn{2}{|c|}{$\textbf{CPOL}$}   \\\cline{2-7}
    	 &	Location	&	Period	&	Location	&	Period	&	Location	&	Period 	\\
    \hline
    School &	North	&	2	&	North	&	2	& North	&	3		\\
 Leisure Centre&  	-	&	-	&	-	&	-	& South	&	0		\\
 Council Offices& 	South	&	0	&	South	&	0	& -	&	-		\\
Recycling Centre &	North	&	0	&	North	&	0&	North	&	0		\\
Start Up Incubator & 	North	&	0	&	South	&	0& 	North	&	4		\\
Healthcare  Centre &	South	&	3	&	South	&	3&	South	&	2		\\
Community Centre&	North	&	1	&	South	&	1	&South	&	1		\\
Social Housing &	-	&	-	&	-	&	-	&-	&	-		\\
     \hline
\end{tabu}
\end{table}

\begin{figure}
\centering
\captionsetup[subfigure]{labelformat=empty}
\subfloat {\includegraphics[width =5 in]{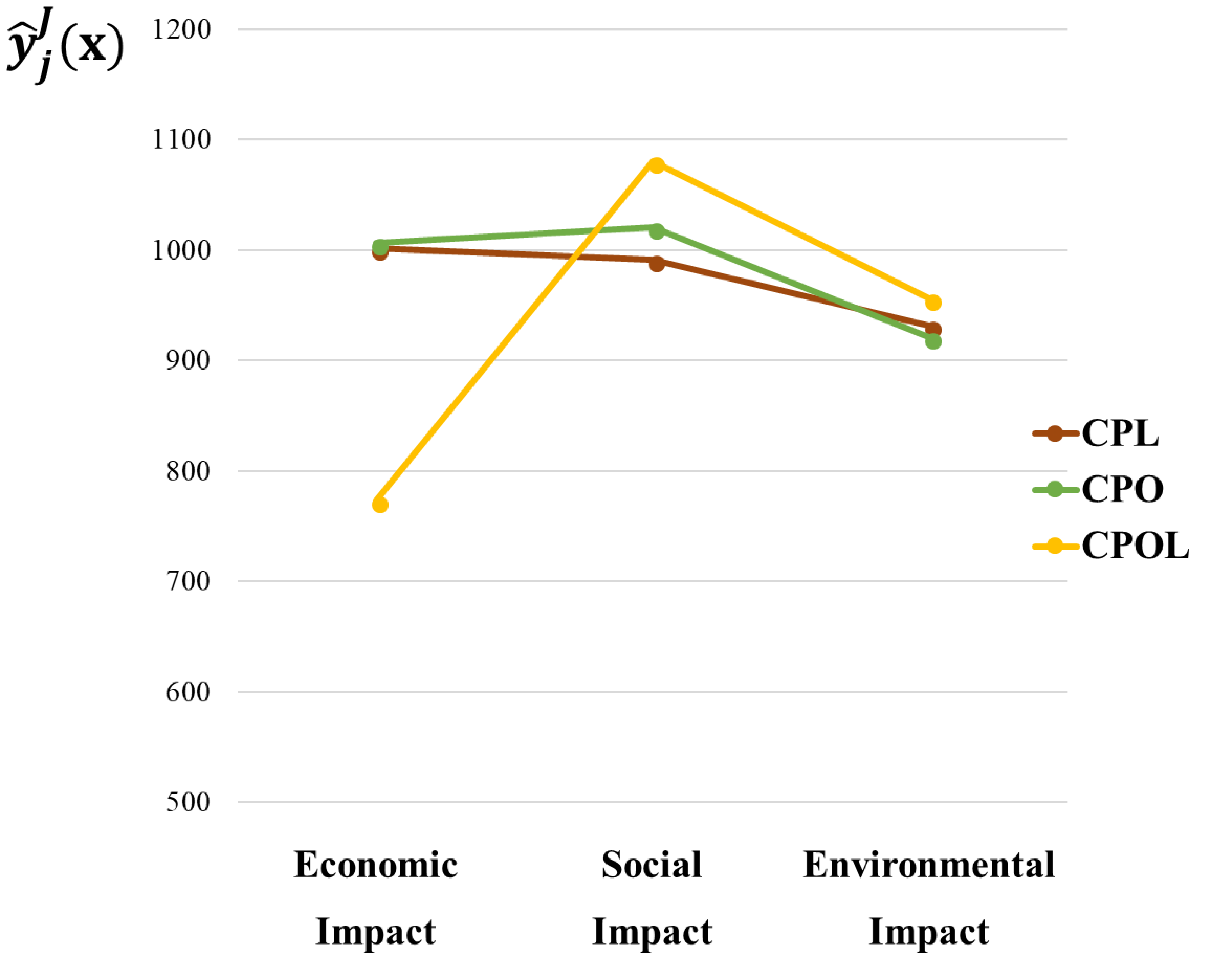}}
\caption{The performances $\widehat{y}^{J}_{j}$ of the best strategies obtained  for each CP approach.}
\label{ObiettiviCP}
\end{figure}

\begin{figure}
\centering
\captionsetup[subfigure]{labelformat=empty}
\subfloat {\includegraphics[width =4.8 in]{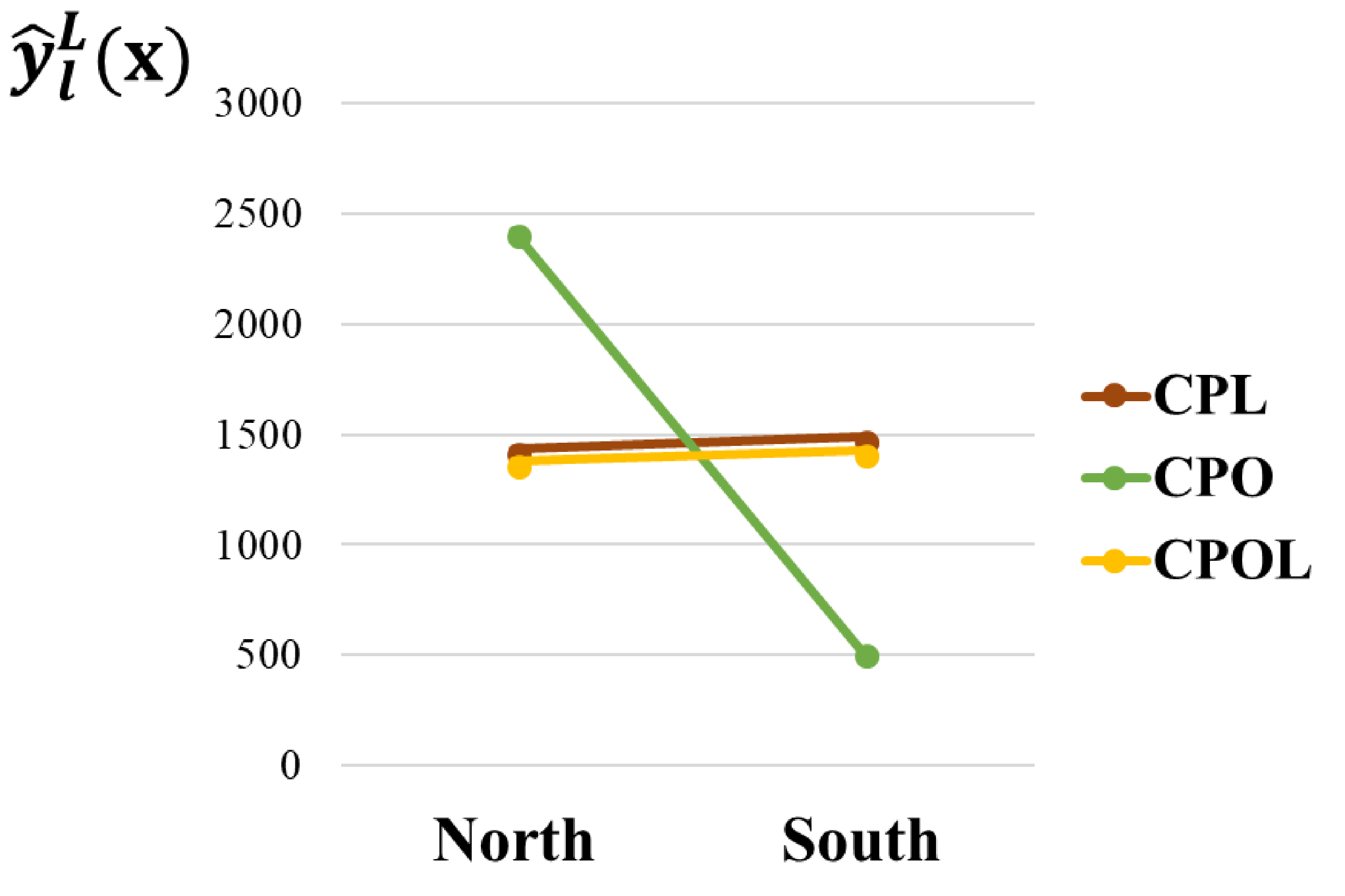}}
\caption{The performances $\widehat{y}^{L}_{l}$ of the best strategies obtained  for each CP approach.}
\label{LocationsCP}
\end{figure}

In Figure \ref{ObiettiviCP} we report the overall performance of the optimal strategy obtained with each of the CP approaches with respect to each criterion, while in Figure \ref{LocationsCP} we show the overall performance with respect to each location. We can see that CPO gives quite balanced values with respect to the  overall performances  $\widehat{y}_j^J$ on considered criteria $j\in J$, while the performances $\widehat{y}_l^L$ with respect to locations $l \in L$ result quite unbalanced. This is because this strategy does not search for a compromise in the values of the differences between the two locations. Nevertheless, CPL and CPOL have similar results. In particular, CPL obtains a solution with a very balanced values $\widehat{y}^{J}_{j}$ between North and South, while CPOL shows a good balance for both criteria and locations.

\subsection{A multiobjective interactive optimization approach}

Interactive Multi-objective Optimization (IMO) methods (for a survey, see~\cite{Miettinenetal}) look for a solution being as much as possible satisfactory for the DM through procedures alternating computation phases (in which multiobjective optimization problems are solved), and dialog phases (in which preference information is collected from the DM). Among the many IMO methods proposed in the literature, we shall take into consideration a method called IMO-DRSA \cite{greco2001}, but of course, any other IMO method can be applied as well. We adopt this method because it collects preference information and gives results of computation in a quite easy and understandable way. The main idea is that the DM is presented wih a list of feasible solutions and she is asked to select one if she is convinced that it is completely satisfactory. In this case the procedure ends. On the contrary, the DM is asked to indicate a set of relatively good solutions in the list, so that a binary partition into classes ``good" and ``others" of the list of proposed solutions is obtained. From such indirect preference information, using the Dominance Based Rough Set Approach (DRSA) \cite{greco2001}, we induce a set of ``if ..., then ...'' decision rules explaining the partition in  ``good" and ``others" in terms of values $g_{j}(\mathbf{x})$ taken for strategy $\mathbf{x}$ by the criteria considered in the multiobjective optimization. More precisely, supposing, without loss of generality, that all criteria are increasing with respect to the preference, the rules are logical statements of the type\\
\begin{center}
if  $g_{j_1}(\mathbf{x})\ge \rho_{j_1}$ and $g_{j_2}(\mathbf{x})\ge \rho_{j_2}$ and $\ldots$ and $g_{j_r}(\mathbf{x})\ge \rho_{j_r}$, then $\mathbf{x}$ is a good solution.
\end{center}
The decision rules so obtained are presented to the DM, that is asked to select she considers the most representative of her preferences. The selected decision rule gives a set of constraints
\begin{center}
$g_{j_1}(\mathbf{x})\ge \rho_{j_1}$, $g_{j_2}(\mathbf{x})\ge \rho_{j_2}$, $\ldots$, $g_{j_r}(\mathbf{x})\ge \rho_{j_r}$
\end{center}
to be added to the current set, so that the solution space is consequently reduced in a region of feasible solutions being more appealing to the DM. From the current set of feasible solutions, another set of representative solutions is built and presented to the DM, so that the cycle starts again, until the DM finds a satisfactory solution.

To apply IMO-DRSA to the specific decision problem represented by our space time model, following \cite{barbati}, the performances $\hat{y}^{sets}_{indices}$ are all transformed in qualitative ordinal evaluations by means of suitable thresholds. With this aim, for each $j \in J$ and each $l \in L$,  the DM is asked to define a set $S_{j,l}$ consisting of $J(h)$ thresholds

$$ S_{j,l}=\{s_{1,j,l};\ldots;s_{J(h),j,l} : s_{1,j,l}<s_{2,j,l}<\ldots<s_{J(h),j,l}\},$$

\noindent permitting to define a set ${\mathcal C}_{h}$ consisting of  $J(h)+1$ qualitative satisfaction classes $C_{a,j,l}$

$$ {\mathcal C}_{h}=\{C_{1,j,l},\ldots,C_{J(h)+1,j,l}\}$$

\noindent such that the greater $a=1,\ldots,J(h)+1$, the more preferred is the project from class $C_{a,j,l}$.
The facilities $i \in I$ are assigned to satisfaction classes $C_{a,j,l} \in {\mathcal C}_h$ according to the following rule: for all $i \in I$
\begin{itemize}
	\item facility $i$ is assigned to class $C_{1,j,l}$ if $y_{ijl}< s_{1,j,l}$;
	\item facility $i$ is assigned to class $C_{a,j,l}$ with $a=2,\ldots,J(h)$, if $s_{a-1,j,l}\le y_{ijl}< s_{a,j,l}$;
	\item facility $i$ is assigned to class $C_{J(h)+1,j,l}$ if  $s_{J(h),j,l}\le y_{ijl}$.
\end{itemize}

\vspace{20mm}
\noindent \emph{Illustrative Example: definition of qualitative valuations}

The council has defined three satisfaction levels for each criterion (see Table \ref{Portfolios_levels}). We suppose the same levels are considered for all criteria in all the locations. In this way we can define our satisfaction classes: ``weakly satisfactory", ``satisfactory", ``very satisfactory", and ``extremely satisfactory".

 \begin{table}
\centering\caption{Satisfaction levels for the three criteria in both the locations considered in the illustrative example.}\label{Portfolios_levels}
\makebox[\linewidth]{
\begin{tabu}{|l||c|c|c|}
\hline
\ {Satisfaction levels}   &   {$Economic Impact$} &   {$Social Impact$} &   {$Environmental Impact$}  \\
       \hline
           \hline
\textit{$s_1$}: Satisfactory 	&	20	&	20	&	20		\\
\textit{$s_2$}: Very satisfactory	&	35	&	35	&	35		\\
\textit{$s_3$}:	Extremely satisfactory &	55	&	55	&	55		\\
\hline
\end{tabu}
}
\end{table}

For each facility $i \in I$ with respect to each location $l \in L$ and to each criterion $j \in J$ we have:

\begin{itemize}
	\item  facility $i$ is ``weakly satisfactory" if $y_{ijl} < 20$;
	\item  facility $i$  is ``satisfactory" if $  20 \leq y_{ijl}< 35$;
	\item  facility $i$  is ``very satisfactory" if $35 \leq y_{ijl}  < 55$;
	\item facility $i$  is ``extremely satisfactory" if $55 \leq  y_{ijl}$. $\diamond$
\end{itemize}

For each strategy $\mathbf{x}$,  in each location $l \in L$, for each criterion $j \in J$, each satisfaction level $s_{a,j,l}\in S_{j,l}$, we consider the set of facilities attaining threshold $s_{a,j,l}$:

$$ P_{a,j,l}(\mathbf{x}) =\{i \in I: y_{ijl}(\mathbf{x})\geq s_{a,j,l}\}.$$

\noindent In simple words, considering the qualitative scale given in the above example, with respect to criterion $j$, for the strategy $\mathbf{x}$,
\begin{itemize}
	%\item $P_{1,j,l}(\mathbf{x})$  is the set of weakly satisfactory facilities,
	\item $P_{1,j,l}(\mathbf{x})$  is the set of satisfactory facilities,
	\item $P_{2,j,l}(\mathbf{x})$  is the set of very satisfactory facilities,
	\item $P_{3,j,l}(\mathbf{x})$  is the set of extremely satisfactory facilities.
\end{itemize}

\noindent For the sake of simplicity, in what follows we shall write $|P_{a,j,l}(\mathbf{x})|$, as $F_{a,j,l}(\mathbf{x})$.

We can consider the following three main formulations of our space-time multiobjective optimization problem:
\begin{itemize}
	\item a \emph{location-oriented multiobjective optimization} in which the objective functions are the sums on all considered criteria of the number of activated facilities attaining an evaluation of at least level $a, a=1, \ldots, h$ in a given location $l\in L$, that is, in the above example, for each facility $i\in I$:
		\begin{itemize}
		\item the number of facilities at least satisfactory for the first criterion plus the analogous number for the second criterion and so on until the last criterion;
		\item the number of facilities at least very satisfactory for the first criterion plus the analogous number for the second criterion and so on until the last criterion;
		\item the number of facilities extremely satisfactory for the first criterion plus the analogous number for the second criterion and so on until the last criterion.
			\end{itemize}
				Therefore the location oriented multiobjective optimization problem can be formulated as
		$$\max \sum_{j\in J}{F_{a,j,l}}({\mathbf{x}}), \ \forall l \in L , \ \forall s_{a,j,l} \in S_{j,l}$$
\noindent under the constraints (\ref{ConstraintBudget}) and (\ref{ConstraintOnlyOneFac}), and the other possible constraints of the original problem.
\item  a \emph{criterion oriented multiobjective optimization} in which the objective functions are the sums on all considered locations of the number of activated facilities of at least level $a, a=1, \ldots, h$, for a given criterion $j\in J$; that is, in the above example, for each criterion $j \in J$:
\begin{itemize}
		\item the number of facilities at least satisfactory in the first location plus the analogous number in the the second location and so on until the last location;
		\item the number of facilities at least very satisfactory in the first location plus the analogous number in the the second location and so on until the last location;
		\item the number of facilities extremely satisfactory for the first location plus the analogous number for the second location and so on until the last location.
			\end{itemize}
				Therefore the criterion oriented multiobjective optimization problem can be formulated as
		$$\max \sum_{l\in L}{F_{a,j,l}}({\mathbf{x}}), \ \forall j \in J , \ \forall s_{a,j,l} \in S_{j,l}$$,
\noindent under the constraints (\ref{ConstraintBudget}) and (\ref{ConstraintOnlyOneFac}), and the other possible constraints of the original problem.
\item a \emph{criterion and location oriented multiobjective optimization} in which the objective functions are combinations of one location $l \in L$, one criterion $j \in J$ and the number of activated facilities of at least level $a, a=1, \ldots, h$; that is, in the above example, for each criterion $j \in J$ and $l \in L$:
\begin{itemize}
		\item the number of facilities at least satisfactory;
		\item the number of facilities at least very satisfactory;
		\item the number of facilities extremely satisfactory.
			\end{itemize}
				Therefore the criterion and oriented multiobjective optimization problem can be formulated as
		$$\max {F_{a,j,l}}({\mathbf{x}}), \ \forall j \in J, \forall l \in L\ \forall s_{a,j,l} \in S_{j,l}$$
\noindent under the constraints (\ref{ConstraintBudget}) and (\ref{ConstraintOnlyOneFac}), and the other possible constraints of the original problem.

\end{itemize}

\subsection{Illustrative example: application of IMO-DRSA}

Let us apply the IMO-DRSA to the decision problem introduced in Section \ref{example}. Taking into account the evaluations of the projects with respect to considered criteria shown in Table \ref{Facilities_evaluation} and the thresholds in Table \ref{Portfolios_levels}, we get the evaluations in ordinal qualitative terms shown in Table \ref{Qualitative_Facilities_evaluation}, where WS, S, ES and VS are representing  our satisfaction classes ``weakly satisfactory", ``satisfactory", ``very satisfactory", and ``extremely satisfactory", respectively.

\begin{table}
\centering\caption{Qualitative ordinal evaluations on three criteria in each location of the eight facilities considered in the illustrative example.}\label{Qualitative_Facilities_evaluation}

\begin{tabu}{|l||c|c|c|c|c|c|}
\hline
   \multirow{2}{*}{Facilities}&   \multicolumn{2}{|c|}{$Economic Impact$} &  \multicolumn{2}{|c|}{$Social Impact$} &   \multicolumn{2}{|c|}{$Environmental Impact$} \\\cline{2-7}

    	 &	North	&	South	&	North	&	South	&	North       	&	South 	\\
      \hline
School & 	 S&	 S&	ES&	ES&	 S&	 S\\
Leisure Centre& 	VS&	VS&	ES&	ES&	VS&	 S\\
Council Offices& 	WS&	WS&	 S&	 S&	 S&	 S\\
Recycling Centre & 	ES&	ES&	ES&	ES&	ES&	ES\\
Start Up Incubator  &	ES&	ES&	WS&	WS&	WS&	WS\\
Healthcare  Centre  &	WS&	WS&	WS&	WS&	VS&	ES\\
Community Centre& 	 S&	 S&	ES&	VS&	 S&	VS\\
Social Housing & 	WS&	 S&	ES&	ES&	WS&	WS\\
     \hline
\end{tabu}
\end{table}

In a perspective of location oriented multiobjective optimization, each portfolio is evaluated in terms of facilities at least satisfactory, at least very satisfactory and extremely satisfactory in the North and in the South. In the first iteration, the six representative strategies  presented in Table \ref{Portfolios_first_iteration} are  shown to the DM,  where $\mathcal{F}_{a,l}({\mathbf{x}})= \sum_{j\in J}{F_{a,j,l}}({\mathbf{x}})$. For example $\mathcal{F}_{1,1}({\mathbf{x}})= \sum_{j\in J}{F_{1,j,1}}({\mathbf{x}})$ is the number of all the facilities that have a contribution for each criterion at least satisfactory.
In Table \ref{SoluzioniPrimaIte} we report here the corresponding strategies. Let us underline that to facilitate the understanding of the solution for the DMs, each strategy could be presented to the DM, with some graphs representing, with histograms, the values of $\mathcal{F}_{a,l}({\mathbf{x}})$ as shown in \cite{barbati}. For the sake of the space we do not report here these representations.

\begin{table}
\centering \caption{The set of non-dominated strategies presented to the DM in the first iteration.}\label{Portfolios_first_iteration}
\makebox[\linewidth]{
\begin{tabu}{|c||c|c|c|c|c|c|c|}
%\hline
 %& \multicolumn{4}{c|}{\textbf{Criteria}}\\
%\hline
\hline
\rowfont{\scriptsize}
\small{{Strategy} }   &   {$\mathcal{F}_{1,1}$} &   {$\mathcal{F}_{1,2}$} &   {$\mathcal{F}_{2,1}$} & {$\mathcal{F}_{2,2}$}& {$\mathcal{F}_{3,1}$} & {$\mathcal{F}_{3,2}$}&\small{{Class}}\\
       \hline
       \hline
\textit{ST}1	&12&	0&	11&	0&	5&	0&		*	\\
\textit{ST}2	&0&	14&	0&	9&	0&	6&		Good	\\
\textit{ST}3	&12&	0&	11&	0&	5&	0&		*	\\
\textit{ST}4	&0&	13&	0&	10&	0&	6&		Good	\\
\textit{ST}5	&8&	2&	7&	1&	6&	1&		*	\\
\textit{ST}6	&1&	12&	1&	9&	0&	6&		Good	\\
\hline
\end{tabu}
}
\end{table}

\begin{table}
\centering\caption{Position and period of the activated facilities for each strategy.}\label{SoluzioniPrimaIte}

\begin{tabu}{|l||c|c||c|c||c|c||c|c||c|c||c|c|}
\hline
   \multirow{2}{*}{Facilities}&   \multicolumn{2}{|c||}{\textit{ST}1} &  \multicolumn{2}{|c||}{\textit{ST}2} &   \multicolumn{2}{|c|}{\textit{ST}3}&   \multicolumn{2}{|c||}{\textit{ST}4} &  \multicolumn{2}{|c||}{\textit{ST}5} &   \multicolumn{2}{|c|}{\textit{ST}6}   \\\cline{2-13}
    	 &	\tiny{Location}	&	\tiny{Period}	&	\tiny{Location}	&	\tiny{Period}	&	\tiny{Location}	&	\tiny{Period} &	\tiny{Location}	&	\tiny{Period}	&	\tiny{Location}	&	\tiny{Period}	&	\tiny{Location}	&	\tiny{Period} 	\\
    \hline
\tiny{School} &	\tiny{North}	&	3	&	\tiny{South}	&	3 &	\tiny{North}	&	3&	\tiny{South}	&	3	&	\tiny{North}	&	2 &	\tiny{South}	&	3		\\
\tiny{Leisure Centre} &  	\tiny{North}	&	0	&	\tiny{South}	&	0	&\tiny{North}	&	0&	\tiny{South}	&	0	&	-	&	- &	\tiny{South}	&	0		\\
\tiny{Council Offices}& 	-	&	-	&	\tiny{South}	&	2&	-	&	-&	-	&	-	&	\tiny{North}	&	4 &	\tiny{North}	&	4		\\
\tiny{Recycling Centre} &	\tiny{North}	&	1	&	\tiny{South}	&	0&	\tiny{North}	&	1&	\tiny{South}	&	0	&	\tiny{North}	&	0 &	\tiny{South}	&	0		\\
\tiny{Start Up Incubator} & 	\tiny{North}	&	4	&	-	&	-&	\tiny{North}	&	0&	-	&	0	&	\tiny{North}	&	3 &	-	&	0		\\
\tiny{Healthcare  Centre} &	\tiny{North}	&	2	&	-	&	-&	\tiny{North}	&	2	&	\tiny{South}	&	2	&	\tiny{North}	&	0 & \tiny{North}	&	2	\\
\tiny{Community Centre} &	\tiny{North}	&	0	&	\tiny{South}	&	1	&\tiny{North}	&	0&	\tiny{South}	&	1	&	\tiny{South}	&	0 &	\tiny{South}	&	1		\\
\tiny{Social Housing} &	-	&	-	&	-	&	-	&-	&	-	&	-	&	-	&	-	&	- &	-	&	-	\\
     \hline
\end{tabu}
\end{table}

The DM is asked if among the strategies shown to her there is one that she considers as completely satisfactory. Since this is not the case, she was asked to select a set of strategies that can be considered as relatively good. Consequently, she indicated strategies ST2, ST4 and ST6. Applying DRSA to this preference information the following decision rules were induced (among parentheses we provide the strategies supporting the corresponding rule):

\begin{itemize}
	\item[Rule 1.1:]  if $\mathcal{F}_{2,2}(\textbf{x})\ge 9$, then strategy $\textbf{x}$  is ``good", \;\;\;\;   \;\;\;\;  \;\;\;\;                ($ST2, ST4, ST6$)
\item[] (if there are at least 9 projects very satisfactory or better in location South with respect to all criteria, then the portfolio is good);
	\item[Rule 1.2:]  if $\mathcal{F}_{3,2}(\textbf{x})\ge 6$, then strategy $\textbf{x}$  is ``good", \;\;\;\;   \;\;\;\;  \;\;\;\;                ($ST2, ST4$, $ST6$)
\item[] (if there are at least 6 projects extremely satisfactory in location South with respect to all criteria, then the portfolio is good);
\item[Rule 1.3:]  if $\mathcal{F}_{1,2}(\textbf{x})\ge 12$, then strategy $\textbf{x}$  is ``good", \;\;\;\;   \;\;\;\;  \;\;\;\;                ($ST2, ST4, ST6$)
\item[] (if there are at least 12 projects very satisfactory or better in location South with respect to all criteria, then the portfolio is good);
\end{itemize}

\noindent The DM selected Rule 1.3 as the most representative for her current aspirations, and the following  constraint  was added to the original optimization problem ${F_{1,2}}({\mathbf{x}})\ge 12$, that is,.
$$\sum_{j\in J}{F_{1,j,2}}({\mathbf{x}})\ge 12$$.

Then, the second sample of weakly non-dominated strategies (shown in Table \ref{Portfolios_second_iteration}) was generated and presented to the DM. For the sake of the space we do not report the correspondent strategies.

  \begin{table}
\centering \caption{The set of non-dominated strategies presented to the DM in the second iteration.}\label{Portfolios_second_iteration}
\begin{tabu}{|c||c|c|c|c|c|c|c|c|}
\hline
\rowfont{\scriptsize}
\small{{Strategy} }   &   {$\mathcal{F}_{1,1}$} &   {$\mathcal{F}_{1,2}$} &   {$\mathcal{F}_{2,1}$} & {$\mathcal{F}_{2,2}$}& {$\mathcal{F}_{3,1}$} & {$\mathcal{F}_{3,2}$}&\small{{Class}}\\
       \hline
       \hline
\textit{ST}1$^{\prime}$&    2&	12&	1&	7&	0&	3&				*	\\
\textit{ST}2$^{\prime}$&	0&	14&	0&	9&	0&	6&				Good	\\
\textit{ST}3$^{\prime}$&	2&	12&	2&	6&	1&	2&				*	\\
\textit{ST}4$^{\prime}$&	0&	13&	0&	10&	0&	6&			Good	\\
\textit{ST}5$^{\prime}$&    1&	12&	1&	8&	1&	5&				*	\\
\textit{ST}6$^{\prime}$&    1&	12&	1&	9&	0&	6&				Good\\
\hline
\end{tabu}
\end{table}

Again, the DM is asked if among the strategies shown to her there is one that she considers as completely satisfactory. Since this is not the case, she was asked to select a set of strategies that can be considered as relatively good. She indicated the strategies apart from ST2$^{\prime}$, ST4$^{\prime}$, ST6$^{\prime}$. Applying DRSA to this preference information the following decision rules were induced (among parentheses we provide the strategies supporting the corresponding rule):

\begin{itemize}
	\item[Rule 2.1:]  if $\mathcal{F}_{2,2}(\textbf{x})\ge 9$ , then strategy $\textbf{x}$  is ``good", \;\;\;\;   \;\;\;\;  \;\;\;\;                (ST$2^{\prime}$, ST$^{\prime}$4, ST$^{\prime}$6)
\item[] (if there are at least 9 projects very satisfactory or better in location South with  respect to all criteria, then the portfolio is good);
\item[Rule 2.2:]  if $\mathcal{F}_{3,2}(\textbf{x})\ge 6$, then strategy $\textbf{x}$  is ``good", \;\;\;\;   \;\;\;\;  \;\;\;\;                (ST$2^{\prime}$, ST$^{\prime}$4, ST$6^{\prime}$)
\item[] (if there are at least 6 projects extremely satisfactory in location South with respect to all criteria, then the portfolio is good);
\item[Rule 2.3:]  if $\mathcal{F}_{1,2}(\textbf{x})\ge 13$, then strategy $\textbf{x}$  is ``good", \;\;\;\;   \;\;\;\;  \;\;\;\;                (ST$2^{\prime}$, ST$4^{\prime}$)
\item[] (if there are at least 13 projects very satisfactory in location South with respect to all criteria, then the portfolio is good);
\end{itemize}

\begin{table}
\centering \caption{A set of non-dominated strategies presented to the DM in the third iteration.}\label{Portfolios_third_iteration}
\makebox[\linewidth]{
\begin{tabu}{|c||c|c|c|c|c|c|c|}
%\hline
 %& \multicolumn{4}{c|}{\textbf{Criteria}}\\
%\hline
\hline
\rowfont{\scriptsize}
\small{{Strategy} }   &   {$\mathcal{F}_{1,1}$} &   {$\mathcal{F}_{1,2}$} &   {$\mathcal{F}_{2,1}$} & {$\mathcal{F}_{2,2}$}& {$\mathcal{F}_{3,1}$} & {$\mathcal{F}_{3,2}$}&\small{{Class}}\\
       \hline
       \hline
\textit{ST}1$^{\prime\prime}$&   1&	12&	1&	9&	0&	6&			Good	\\
\textit{ST}2$^{\prime\prime}$&	0&	14&	0&	9&	0&	6&			*	\\
\textit{ST}3$^{\prime\prime}$&	1&	12&	1&	9&	0&	6&		*	\\
\textit{ST}4$^{\prime\prime}$&	0&	13&	0&	10&	0&	6&			*	\\
\textit{ST}5$^{\prime\prime}$&   0&	14&	0&	9&	0&	6&		*	\\
\textit{ST}6$^{\prime\prime}$ &  0&	13&	0&	10&	0&	6&		*\\
\hline
\end{tabu}
}
\end{table}

\noindent The DM selected Rule 2.1 as the most representative for her current aspirations, and thus, the following constraint was added to the original optimization problem and to the constraints added to the previous interaction ${F_{2,2}}({\mathbf{x}})\ge 9$, that is,
$$\sum_{j\in J}{F_{2,j,2}}({\mathbf{x}})\ge 9$$.

Then, the third sample of weakly non-dominated strategies shown in Table \ref{Portfolios_third_iteration} was generated and presented to the DM.

At this point the DM declares to be satisied by the strategy ST1$^{\prime\prime}$ and the procedure stops.
%%%%%%%%%%%%%%%%%%%%%%%%%%%%%%%%%%%%%%%%%
\section{Uncertainty and plurality of stakeholders}\label{IncertezzaPlur}
%%%%%%%%%%%%%%%%%%%%%%%%%%%%%%%%%%%%%%%%%%%%%%%%%%%%
Two features that affect many real world problems are related to the uncertainty of the performances expected from activation of facilities \cite{Mild,Vilkkumaa} and to the presence of a plurality of stakeholders \cite{SaloHAm}. In the following we introduce these further elements in our model.

\subsection{Uncertainty}
We model the uncertainty related to the performances of facilities from  $I$ with respect to criteria $j$ from $J$ taking into account a set of states of nature related to the period $t \in T$ and to the states of nature realized in previous periods. Therefore we denote by
$$s_{(t,h_1,\ldots,h_t)}, t \in T-{0}$$
a state of nature taking place in period $t$ in the sequence of previous states of nature,
$$s_{(1,h_1)},s_{(2,h_1,h_2)},\ldots,s_{(t-1,h_1,h_2,\ldots,h_{t-1})}.$$
For all $t \in T - \{0\}$ and for all path $s_{(1,h_1)},s_{(2,h_1,h_2)},\ldots, s_{(t-1,h_1,h_2,\ldots,h_{t-1})}$, let us denote by
$$p^C(s_{(t,h_1,h_2,\ldots,h_{t})})$$
the probability of $s_{(t,h_1,h_2,\ldots,h_{t})}$ conditioned to the path of previous states of nature
$$s_{(2,h_1,h_2)}, \ldots, s_{(t-1,h_1,h_2,\ldots,h_{t-1})}.$$
In other words, $p^C(s_{(t,h_1,h_2,\ldots,h_{t})})$ is the probability of realization in period $t$ of $s_{(t,h_1,h_2,\ldots,h_{t})}$ if, in period $t-1$ state of nature $s_{(t,h_1,h_2,\ldots,h_{t-1})}$, is realized.
Consequently, the (non conditioned) probability of the state of nature $s_{(t,h_1,\ldots,h_t)}$ is given by
$$p(s_{(t,h_1,\ldots,h_t)})=p^C(s_{(1,h_1)})\times p^C(s_{(2,h_1,h_2)})\times \ldots \times p^C(s_{(t,h_1,\ldots,h_t)}).$$

\begin{figure}
\centering
\captionsetup[subfigure]{labelformat=empty}
\subfloat {\includegraphics[width = 6in]{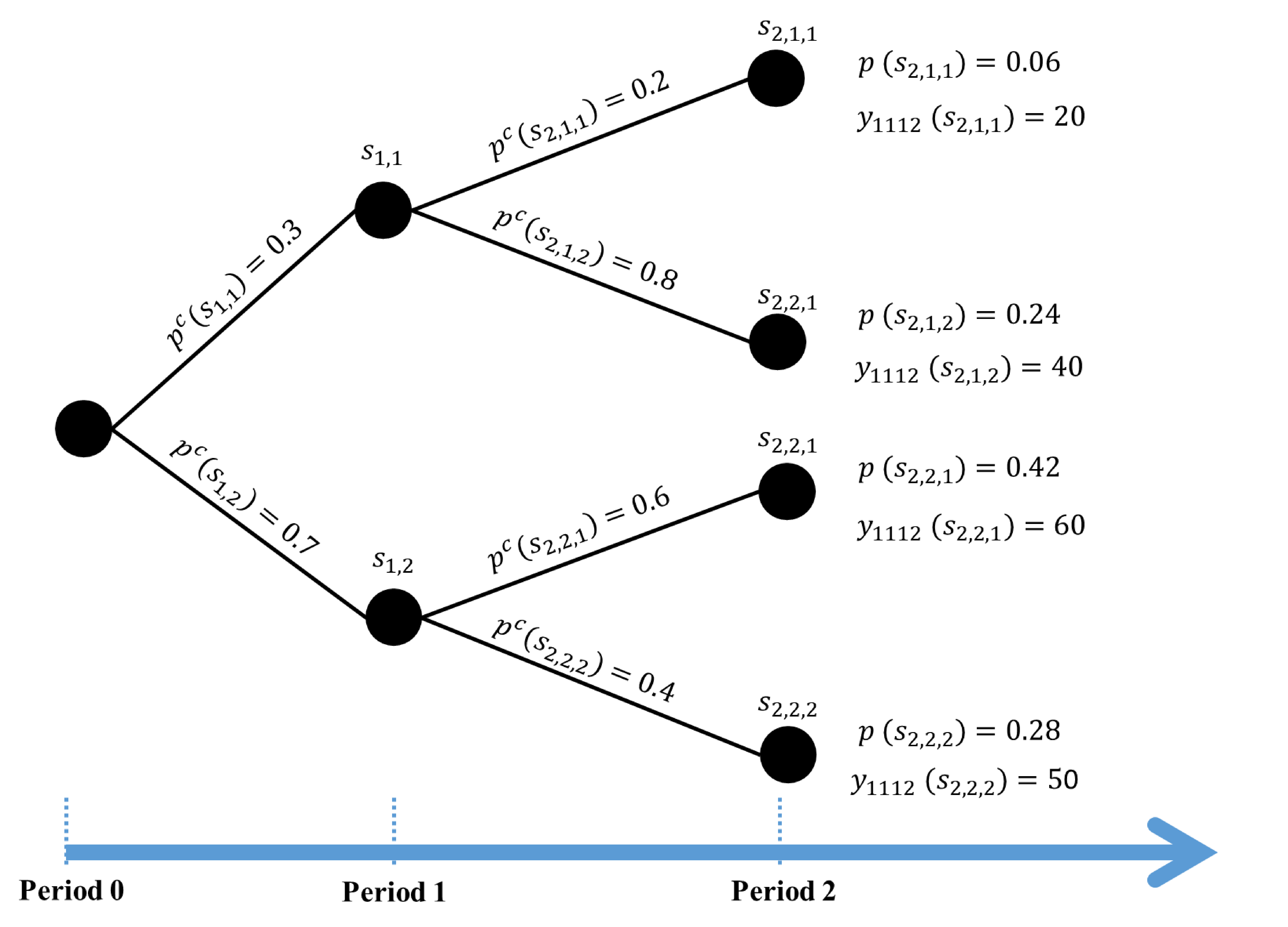}}
\caption{Example of a probability distribution of the performances.}
\label{DiagramTree}
\end{figure}

For instance, let us consider the example in Figure \ref{DiagramTree} regarding the first facility, in the first location and for the first period. We have 2 periods and, for each period, we have two possible states of nature. Every state of nature is associated to a node of the diagram tree and the probability of each state of nature $p^C(s_{(t,h_1,h_2,\ldots,h_{t})})$ is reported on the arc entering each node. 

In this context we denote by $y_{ijlt}(s_{(t,h_1,\ldots,h_t)})$ the performance of facility $i \in I$ with respect to criterion $j \in J$ in location $l \in L$ at time $t \in T$ if the state of nature $s_{(t,h_1,\ldots,h_t)}$ is realized.

Taking into account the probabilities $p(s_{(t,h_1,\ldots,h_t)})$, we can compute the expected value of the performance of facility $i \in I$ with respect to criterion $j \in J$ in location $l \in L$ at time $t\in T-\{0\}$ as follows:
$$E_p(y_{ijlt})=\sum_{s_{(t,h_1,h_2,\ldots,h_{t})}\in S_t} y_{ijlt}(s_{(t,h_1,\ldots,h_t)}) \times p(s_{(t,h_1,\ldots,h_t)}),$$
where $S_t$ denotes the set of possible states of nature in period $t$.

In our didactic example, in Figure \ref{DiagramTree}, the expected value $E_p(y_{ijl2})$ can be calculated as follows:
$$E_p(y_{ijl2})=0.06\times 20 + 0.24\times 40 + 0.42 \times 60 + 0.28 \times 50 = 50$$,.

Given a strategy $\mathbf{x}$, the expected value $E_p(y^{IJLT}_{ijlt}(\mathbf{x}))$ of the performance of criterion $j \in J$ in period $t\in T-\{0\}$ from facility $i \in I$  is given by
$$E_p(y^{IJLT}_{ijlt}(\mathbf{x}))=\sum_{\tau=0}^{t-1} x_{il\tau}E_p(y_{ijlt}).$$
Analogously, the expected value of the discounted performance $y^{IJLT}_{ijlt}(\mathbf{x})$ is the following
$$E_p(\widehat{y}^{IJLT}_{ijlt}(\mathbf{x}))= E_p(y^{IJLT}_{ijlt}(\mathbf{x}))v(t)=\sum_{\tau=0}^{t-1} x_{il\tau}E_p(y_{ijlt})v(t).$$

Moreover, the expected value of all the other interesting values obtained from the values $y^{IJLT}_{ijlt}(\mathbf{x})$, can be easily obtained using $E_p(y^{IJLT}_{ijlt}(\mathbf{x}))$ instead of $y^{IJLT}_{ijlt}(\mathbf{x})$, as well as the corresponding discounted values can be obtained using $E_p(\widehat{y}^{IJLT}_{ijlt}(\mathbf{x}))$ instead of $\widehat{y}^{IJLT}_{ijlt}(\mathbf{x})$. For example the expected value of the
the global performance of the strategy $\mathbf{x}$  with respect to criterion $j \in J$ in location $l \in L$ at time $t\in T-\{0\}$ is
$$E_p(y^{JLT}_{jlt}(\mathbf{x}))=\sum_{i \in I}E_p(y^{IJLT}_{ijlt}(\mathbf{x})) =\sum_{i \in I} \sum_{\tau=0}^{t-1} x_{il\tau}E_p(y_{ijlt})$$
and its discounted value is
$$E_p(\widehat{y}^{JLT}_{jlt}(\mathbf{x}))=\sum_{i \in I}E_p(\widehat{y}^{IJLT}_{ijlt}(\mathbf{x})) =\sum_{i \in I} \sum_{\tau=0}^{t-1} x_{il\tau}E_p(y_{ijlt})v(t).$$

In first approximation, the problem to be handled is to select the strategy $\mathbf{x}$ maximizing the expected value of the discounted overall performance of strategy $\mathbf{x}$ taking into account all facilities $i \in I$, all criteria $j \in J$, all locations $l \in L$ and all periods $t\in T-\{0\}$ that is
$$E_p(\widehat{y}(\mathbf{x}))=\sum_{i \in I}\sum_{j \in J}\sum_{l \in L}\sum_{t \in T-\{0\}}\sum_{\tau=0}^{t-1} w_j x_{il\tau}E_p(y_{ijlt})v(t).$$

Of course, also in this case one can handle the selection of the most preferred strategy by defining some compromise programming problem analogous to those ones illustrated in Section \ref{multiobjective}. One can use also some  interactive multiobjective optimization method, such as the IMO DRSA again introduced in Section 4. In this perspective, to deal with uncertainty performances and time preferences using with DRSA, one can follow the approach proposed in  \cite{Greco2010}.

\vspace{5mm}
\noindent \emph{Illustrative Example: Uncertainty}

\begin{figure}
\centering
\captionsetup[subfigure]{labelformat=empty}
\subfloat {\includegraphics[width = 6in]{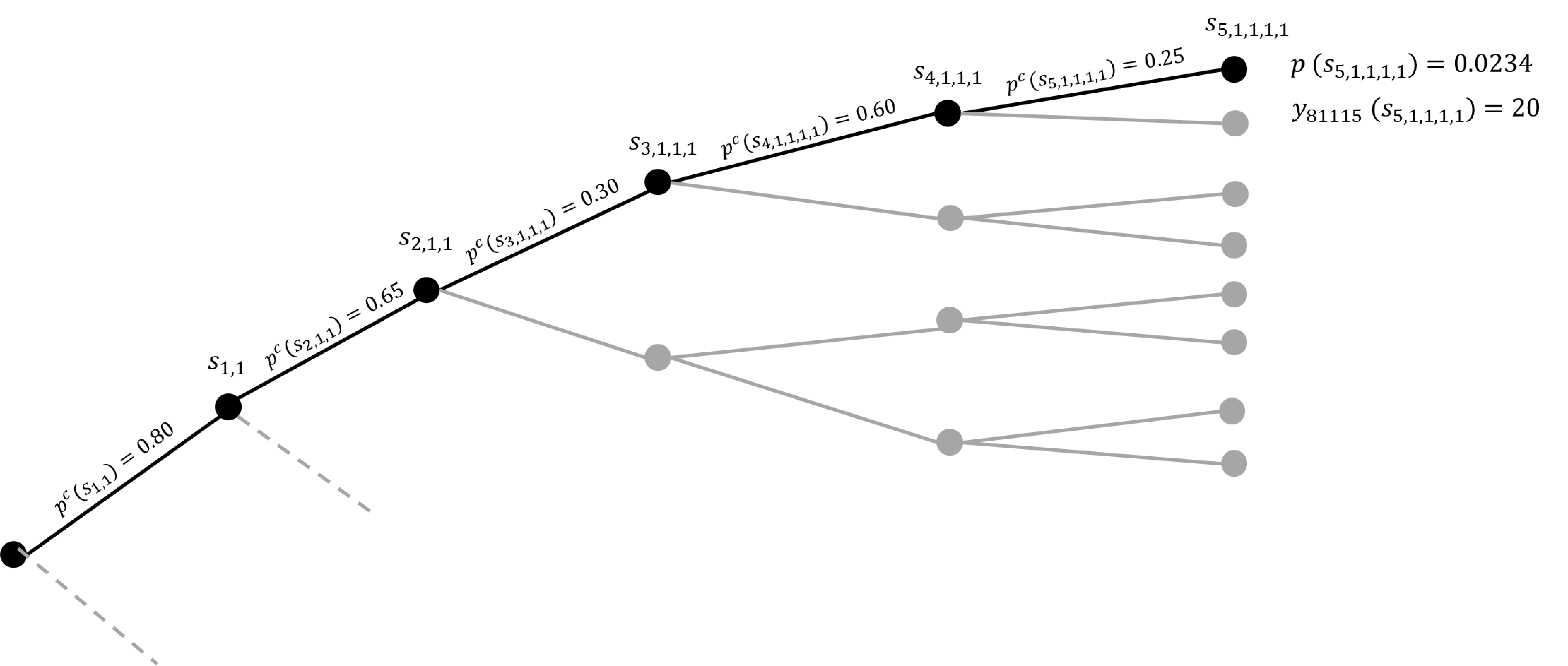}}
\caption{Probability distribution of the performances for the facility Social Housing with respect to the economic aspects.}
\label{EsempioTree}
\end{figure}

In order to show how uncertainty can be taken into account with the proposed approach, we reconsider, for example, the performances  related to  the facility \emph{Social Housing}, with respect to  the economic criterion, having a probability distribution in both the locations with two possible alternative states of nature in each period.  For the sake of space limit, we reported in Figure \ref{EsempioTree} only one branch of the tree. Following the path highlighted in bold black, we can compute   $p(s_{(5,1,1,1,1,1)})=0.80\times 0.65\times 0.30 \times 0.60 \times 0.25=0.023$ to which is associate the performance $y_{811}(s_{(5,1,1,1,1,1)})=76$. The other evaluations and the other probabilities are listed in Table \ref{ProbEval}.

\begin{figure}
\centering
\captionsetup[subfigure]{labelformat=empty}
\subfloat {\includegraphics[width = 6in]{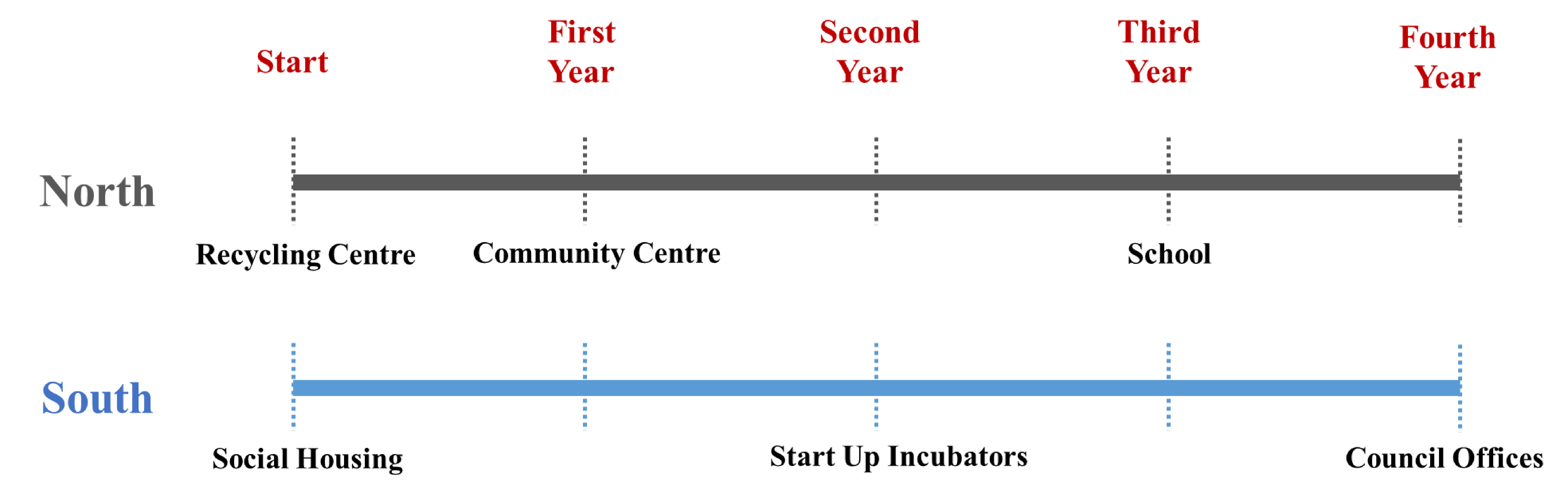}}
\caption{Optimal solution obtained by maximizing the expected value of the discounted overall performance $E_p(\hat{y}(\mathbf{x}))$.}
\label{ProbabilitySolution}
\end{figure}

The evaluations $y_{ijl}$ of our illustrative example in Table \ref{Facilities_evaluation} remain the same, apart from $y_{811}$ and $y_{812}$ changed in $E_p(y_{811})=73$ and $E_p(y_{812})=76$, respectively. Maximizing   the expected value of the discounted overall performance $E_p(\hat{y}(\mathbf{x}))$ we obtain the most preferred solution shown in Figure \ref{ProbabilitySolution}. We can note that the facility Social Housing has to be activated in the first period; in fact, its economic evaluation is much improved in comparison to the not probabilistic scenario and this determines its entrance in the optimal strategy; nevertheless some very low evaluations of the facilities are taken into account also in our probabilistic scenario.

\begin{table}
\centering \caption{Performances $y_{811}(s_{(t,h_1,\ldots,h_t)})$  and corresponding probabilities $y_{812}(s_{(t,h_1,\ldots,h_t)})$  for each possible state of nature in the final period .}\label{ProbEval}
\begin{tabu}{|l||c|c|c|c|}
\hline
\cline{2-5}
\multirow{2}{*}{  State of Nature} &   \multicolumn{2}{|c|}{$y_{811}$} &  \multicolumn{2}{|c|}{$y_{812}$}   \\
\cline{2-5}
& Probabilities & Performances& Probabilities & Performances\\
\hline
$s_{(5,1,1,1,1,1)}$		&	0.0234	&	76	&	0.0234	&	89	\\
$s_{(5,1,1,1,1,2)}$		&	0.0702	&	64	&	0.0702	&	92	\\
$s_{(5,1,1,1,2,1)}$		&	0.01872	&	96	&	0.01872	&	84	\\
$s_{(5,1,1,1,2,2)}$		&	0.04368	&	78	&	0.04368	&	95	\\
$s_{(5,1,1,2,1,1)}$		&	0.0364	&	81	&	0.0364	&	78	\\
$s_{(5,1,1,2,1,2)}$		&	0.1456	&	86	&	0.1456	&	93	\\
$s_{(5,1,1,2,2,1)}$		&	0.1092	&	66	&	0.1092	&	17	\\
$s_{(5,1,1,2,2,2)}$		&	0.0728	&	69	&	0.0728	&	99	\\
$s_{(5,1,2,1,1,1)}$		&	0.00819	&	78	&	0.00819	&	96	\\
$s_{(5,1,2,1,1,2)}$		&	0.04641	&	64	&	0.04641	&	88	\\
$s_{(5,1,2,1,2,1)}$		&	0.01764	&	81	&	0.01764	&	12	\\
$s_{(5,1,2,1,2,2)}$		&	0.01176	&	67	&	0.01176	&	78	\\
$s_{(5,1,2,2,1,1)}$		&	0.0196	&	90	&	0.0196	&	69	\\
$s_{(5,1,2,2,1,2)}$		&	0.0784	&	81	&	0.0784	&	87	\\
$s_{(5,1,2,2,2,1)}$		&	0.049	&	67	&	0.049	&	79	\\
$s_{(5,1,2,2,2,2)}$		&	0.049	&	95	&	0.049	&	94	\\
$s_{(5,2,1,1,1,1)}$		&	0.00612	&	39	&	0.00612	&	15	\\
$s_{(5,2,1,1,1,2)}$		&	0.02448	&	68	&	0.02448	&	92	\\
$s_{(5,2,1,2,2,1)}$		&	0.00162	&	76	&	0.00162	&	15	\\
$s_{(5,2,1,2,2,2)}$		&	0.00378	&	26	&	0.00378	&	77	\\
$s_{(5,2,1,1,1,1)}$		&	0.0084	&	80	&	0.0084	&	69	\\
$s_{(5,2,1,1,1,2)}$		&	0.0336	&	70	&	0.0336	&	67	\\
$s_{(5,2,1,2,2,1)}$		&	0.0252	&	94	&	0.0252	&	93	\\
$s_{(5,2,1,2,2,2)}$		&	0.0168	&	43	&	0.0168	&	12	\\
$s_{(5,2,2,1,1,1)}$		&	0.00432	&	62	&	0.00432	&	75	\\
$s_{(5,2,2,1,1,2)}$		&	0.01008	&	44	&	0.01008	&	88	\\
$s_{(5,2,2,2,2,1)}$		&	0.00384	&	65	&	0.00384	&	77	\\
$s_{(5,2,2,2,2,2)}$		&	0.00576	&	26	&	0.00576	&	10	\\
$s_{(5,2,2,1,1,1)}$		&	0.00448	&	66	&	0.00448	&	48	\\
$s_{(5,2,2,1,1,2)}$		&	0.01792	&	51	&	0.01792	&	87	\\
$s_{(5,2,2,2,2,1)}$		&	0.01176	&	42	&	0.01176	&	15	\\
$s_{(5,2,2,2,2,2)}$		&	0.02184	&	41	&	0.02184	&	96	\\
\hline
\end{tabu}
\end{table}

\subsection{Plurality of stakeholders}
In planning problems we usual  have a plurality of stakeholders such as municipality, building companies, association of citizens, trade union and so on \cite{Montibeller}. Therefore it is reasonable to generalize our model to the presence of different perspectives and preferences expressed by different stakeholders. Here we present a basic approach of group decisions to our space-time model. Of course, more complex  approaches can be considered.  The basic idea is to assume a different weights vector for each stakeholder. Let us suppose that we have $K=\{1, \ldots, k, \ldots, b\}$ stakeholders. We consider weights $w_{jk} \ge 0$, with $w_{1k}+\ldots+w_{qk}=1$,  where $w_{jk}$ represents the weight assigned to criterion $j$ from stakeholder $k$.   We  also introduce a central planner that defines a compromise solution giving a weight $z_k \ge 0$ representing the importance of each stakeholder, with $z_{1}+\ldots+z_{b}=1$.

In this way, among the great plurality of performances $y^{sets}_{indices}(\mathbf{x})$ defined in Section \ref{proposed model} we can reformulate some of them and add others as follows:

\begin{itemize}

\item the overall performance of facility $i \in I$ in location $l \in L$ in period $t\in T-\{0\}$ taking into account all criteria, for  stakeholder $k \in K$ is
$$y^{ILTK}_{iltk}(\mathbf{x})=\sum_{j \in J}\sum_{\tau=0}^{t-1} w_{jk} x_{il\tau}y_{ijl},$$

\item the global performance of the strategy $\mathbf{x}$  in location $l \in L$ in period $t\in T-\{0\}$ for stakeholder $k \in K$ is
$$y^{LTK}_{ltk}(\mathbf{x})=\sum_{i \in I} \sum_{j \in J} \sum_{\tau=0}^{t-1} w_{jk} x_{il\tau}y_{ijl}$$

\item the overall performance of facility $i \in I$ in location $l \in L$ in period $t\in T-\{0\}$ taking into account all criteria and all stakeholders is
$$y^{ILT}_{ilt}(\mathbf{x})=\sum_{k \in K}z_k\sum_{j \in J}\sum_{\tau=0}^{t-1} w_{jk} x_{il\tau}y_{ijl},$$

\item the performance of facility $i \in I$  in period $t\in T-\{0\}$ considering all criteria $j \in J$ and all locations $l \in L$ for stakeholder $k \in K$ is
$$y^{ITK}_{itk}(\mathbf{x})=\sum_{j \in J}\sum_{l \in L}\sum_{\tau=0}^{t-1} w_{jk} x_{il\tau}y_{ijl},$$

\item the performance of facility $i \in I$  in location $l \in L$ for stakeholder $k \in K$ considering all criteria $j \in J$ and all periods $t\in T-\{0\}$ is
$$y^{ILK}_{ilk}(\mathbf{x})=\sum_{j \in J}\sum_{t \in T-\{0\}}\sum_{\tau=0}^{t-1} w_{jk} x_{il\tau}y_{ijl},$$

\item the overall performance of strategy $\mathbf{x}$ in period $t\in T-\{0\}$ for stakeholder $k \in K$ considering all facilities $i \in I$, all criteria $j \in J$ and all  locations $l \in L$ is
$$y^{TK}_{tk}(\mathbf{x})=\sum_{i \in I}\sum_{j \in J}\sum_{l \in L}\sum_{\tau=0}^{t-1} w_{jk} x_{il\tau}y_{ijl},$$

\item the overall performances of strategy $\mathbf{x}$  in location $l \in L$ for stakeholder $k \in K$ considering all criteria $j \in J$ and all periods $t\in T-\{0\}$ is
$$y^{LK}_{lk}(\mathbf{x})=\sum_{i \in I}\sum_{j \in J}\sum_{t \in T-\{0\}}\sum_{\tau=0}^{t-1} w_{jk} x_{il\tau}y_{ijl},$$

\item the global performance of the strategy $\mathbf{x}$ for all the stakeholders in location $l \in L$ at time $t\in T-\{0\}$  is
$$y^{LT}_{lt}(\mathbf{x})=\sum_{k \in K}z_k \sum_{i \in I} \sum_{j \in J} \sum_{\tau=0}^{t-1} w_{jk} x_{il\tau}y_{ijl}$$

\item the performance of facility $i \in I$  in period $t\in T-\{0\}$ considering all criteria $j \in J$ and all locations $l \in L$ for all the stakeholders is
$$y^{IT}_{it}(\mathbf{x})=\sum_{k \in K}z_k\sum_{j \in J}\sum_{l \in L}\sum_{\tau=0}^{t-1} w_{jk} x_{il\tau}y_{ijl},$$

\item the performance of facility $i \in I$  in location $l \in L$ considering all criteria $j \in J$ and all periods $t\in T-\{0\}$ for all stakeholders $k\in K$ is
$$y^{IL}_{il}(\mathbf{x})=\sum_{k \in K}z_k \sum_{j \in J}\sum_{t \in T-\{0\}}\sum_{\tau=0}^{t-1} w_{jk} x_{il\tau}y_{ijl},$$

\item the performance of facility $i \in I$ and for stakeholder $k \in K$ with respect to all criteria $j \in J$, all location $l \in L$ and all periods $t\in T-\{0\}$ is
$$y^{IK}_{ik}(\mathbf{x})=\sum_{j \in J}\sum_{l \in L}\sum_{t \in T-\{0\}}\sum_{\tau=0}^{t-1} w_{jk} x_{il\tau}y_{ijl},$$

\item the overall performances of strategy $\mathbf{x}$  in location $l \in L$ considering all criteria $j \in J$, all periods $t\in T-\{0\}$ and all  stakeholders $k\in K$ is
$$y^{L}_{l}(\mathbf{x})=\sum_{k \in K}z_k\sum_{i \in I}\sum_{j \in J}\sum_{t \in T-\{0\}}\sum_{\tau=0}^{t-1} w_{jk} x_{il\tau}y_{ijl},$$

\item the overall performance of strategy $\mathbf{x}$ for stakeholder $k \in K$ taking into account all facilities $i \in I$, all criteria $j \in J$, all locations $l \in L$ and all periods $t\in T-\{0\}$ is
$$y^{K}_{k}(\mathbf{x})=\sum_{i \in I}\sum_{j \in J}\sum_{l \in L}\sum_{t \in T-\{0\}}\sum_{\tau=0}^{t-1} w_{jk} x_{il\tau}y_{ijl}.$$

\item the performance of facility $i \in I$  with respect to all criteria $j \in J$, all location $l \in L$, all periods $t\in T-\{0\}$ and all  stakeholders $k\in K$is
$$y^{I}_{i}(\mathbf{x})=\sum_{k \in K}z_k \sum_{j \in J}\sum_{l \in L}\sum_{t \in T-\{0\}}\sum_{\tau=0}^{t-1} w_{jk} x_{il\tau}y_{ijl},$$

\item the overall performance in period $t\in T-\{0\}$ considering all facilities $i \in I$, all criteria $j \in J$, all  locations $l \in L$ and all  stakeholders $k\in K$  is
$$y^{T}_{t}(\mathbf{x})=\sum_{k \in K}z_k \sum_{i \in I}\sum_{j \in J}\sum_{l \in L}\sum_{\tau=0}^{t-1} w_{jk} x_{il\tau}y_{ijl},$$

\item the overall performance of strategy $\mathbf{x}$ taking into account all facilities $i \in I$, all criteria $j \in J$, all locations $l \in L$ all periods $t\in T-\{0\}$ and all stakeholders is
$$y(\mathbf{x})=\sum_{k \in K}z_k \sum_{i \in I}\sum_{j \in J}\sum_{l \in L}\sum_{t \in T-\{0\}}\sum_{\tau=0}^{t-1} w_{jk} x_{il\tau}y_{ijl}.$$

\end{itemize}

As before, in the first instance, the problem is to define the strategy $\mathbf{x}$ giving the maximum overall discounted performance $\widehat{y}(\mathbf{x})$ subject to the constraints of the problem such as the budget constraints and the activation constraints.  However, the definition of several $\widehat{y}^{sets}_{indices}(\mathbf{x})$ can be an even   richer dashboard that can be handled as of multiobjective optimization of  performances $y^{sets}_{indices}(\mathbf{x})$ and $\widehat{y}^{sets}_{indices}(\mathbf{x})$ for multiple stakeholders.
In addition, to search for the most preferred solution adopting the weighted approach, we can use a Compromise Programming approach dealing with multiple stakeholders (see for i.e., \cite{Phua}).

In our case we characterize our target as the optimal performance $$\widehat{y}^{K^\ast}_{k}=\max_{\mathbf{x}}\widehat{y}^{K}_{k}(\mathbf{x})=\sum_{i \in I}\sum_{j \in J}\sum_{t \in T-\{0\}}\sum_{\tau=0}^{t-1} w_{jk} x_{il\tau}y_{ijl}v(t)$$ that a strategy $\mathbf{x}$  can attain for  stakeholder $k\in K$.
Following \cite{Drezner2}, in order to get a balanced solution, we minimize the maximum deviation $\Delta_k^K$, on the set of stakeholders $k \in K$, defined as
$$
\Delta_k^K (\mathbf{x})=  \frac{\widehat{y}^{K^\ast}_{k} - \widehat{y}^{K}_{k}(\mathbf{x})}{\widehat{y}^{K^\ast}_{k}}.
$$
Then, the distance of the strategy  $ \mathbf{x}$ from the ideal point is $\Delta^K \mathbf{(x)}=\max_{k\in K} \Delta_k^K(\mathbf{x})$. Consequently, $\Delta^k(\mathbf{x})$ is the objective to be minimizied to get the compromise solution searched for.
This compromise optimisation strategy is particularly suitable in case the stakeholders needs some reciprocal concessions between them in order to reach a consensus on a shared decision.

\vspace{5mm}
\noindent \emph{Illustrative Example: Plurality of Stakeholders}

\begin{figure}
\centering
\captionsetup[subfigure]{labelformat=empty}
\subfloat {\includegraphics[width = 6in]{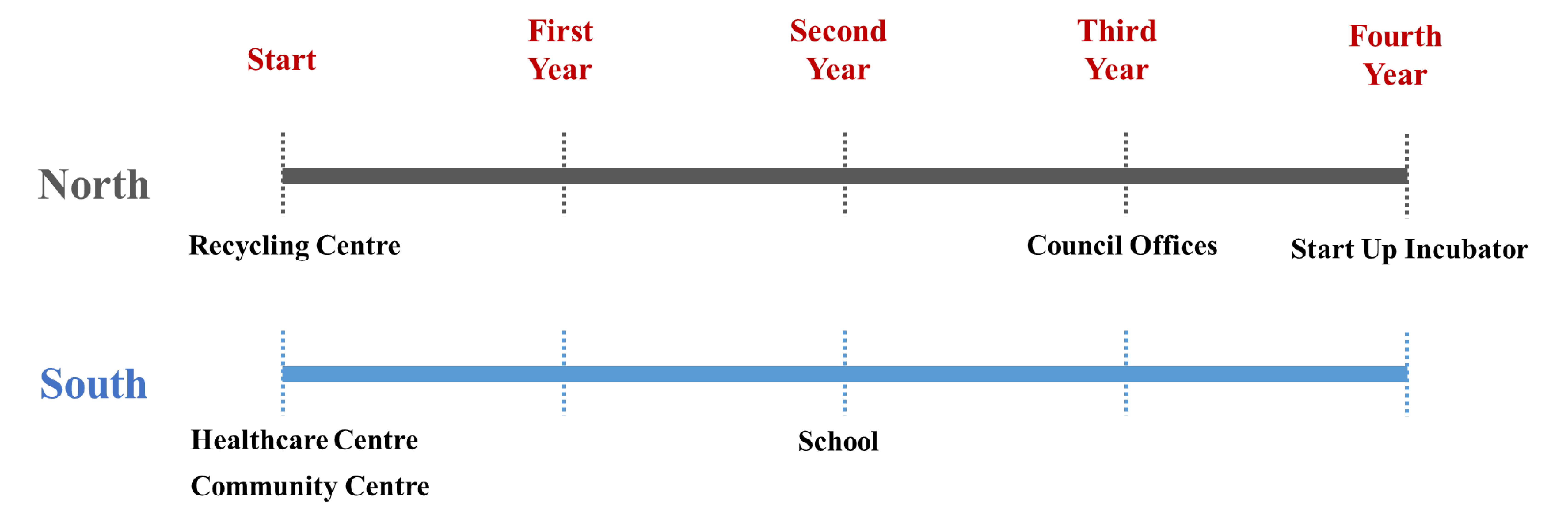}}
\caption{Optimal solution obtained by maximizing the overall performance aggregating the preferences of all the stakeholders with the weights of the central planner.}
\label{MlStkSolution}
\end{figure}

We apply the utility approach to the initial problem defined in Section \ref{example} considering three stakeholders e.g., committees of the council with different interests: Development committee, Planning committee and Government committee. The weights $w_{jk}$ and $z_k$ are reported in Table \ref{weightsmtk}.
Using a utility approach we find the  optimal solution of Figure \ref{MlStkSolution}.

\begin{table}
\centering \caption{Values  for $w_{jk}$ and $z_k$.}\label{weightsmtk}
\begin{tabu}{|l||c|c|c||c|}
\hline
 &	Economic Impact	&	Social Impact	&	Environmental Impact	&	$z_k$	\\
 \hline
Planning committee &  	0.1	&	0.1	&	0.8	&	0.5	\\
 Development committee& 	0.1	&	0.2	&	0.7	&	0.4	\\
Government committee &	0.4	&	0.3	&	0.4	&0.1		\\
     \hline
\end{tabu}
\end{table}

\begin{figure}
\centering
\captionsetup[subfigure]{labelformat=empty}
\subfloat {\includegraphics[width = 7in]{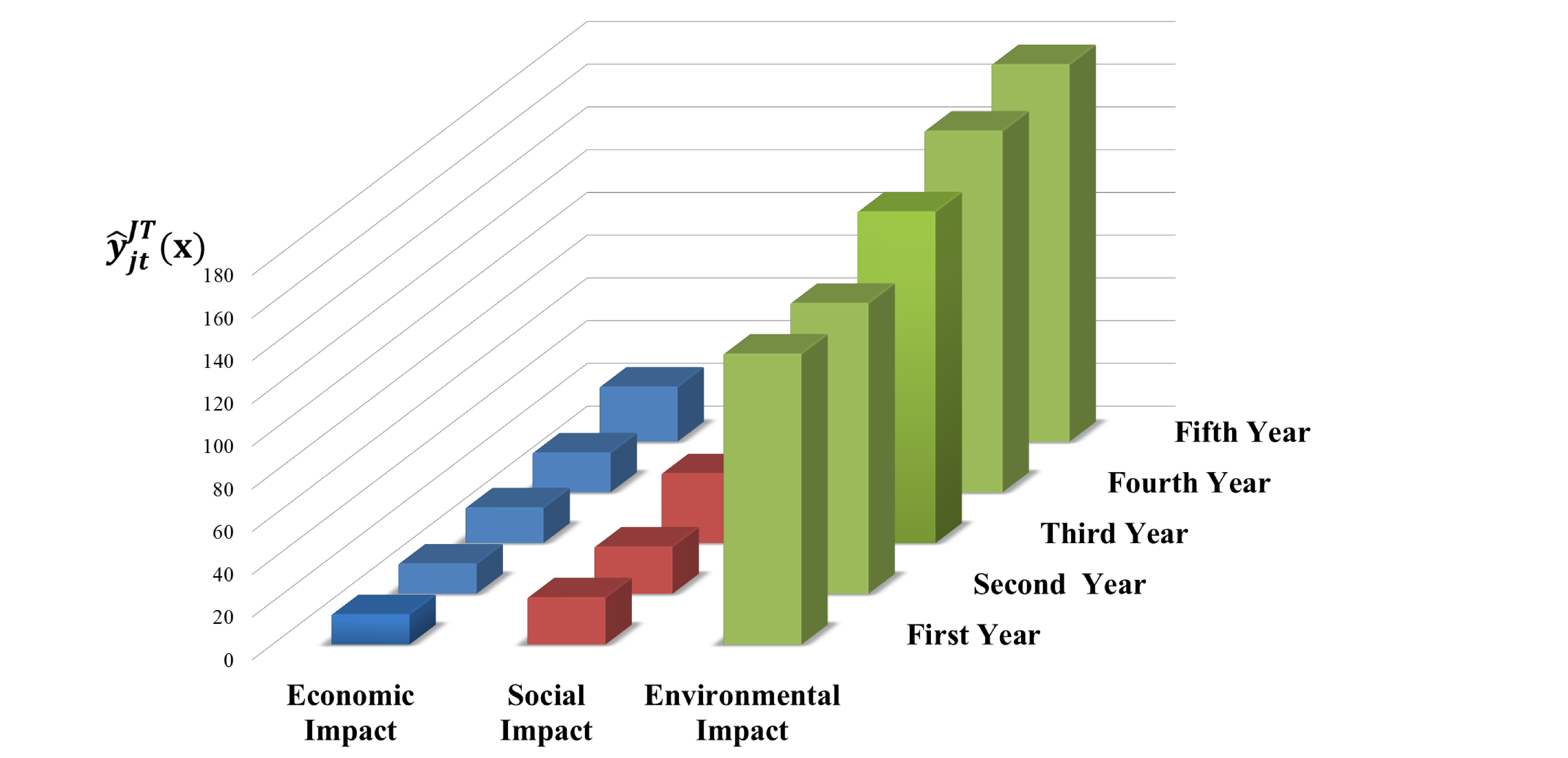}}
\caption{Time distribution of the performances for the optimal  aggregating preferences of all stakeholders, with respect to the considered criteria.}
\label{AttributiTempoMStk}
\end{figure}

From Figure \ref{AttributiTempoMStk} we can see that Environmental Impact is largely more important than the other two criteria because it has a very high weight for the first two committees which are also the most important ones.  This is even more evident if this optimal strategy is compared with the strategy obtained with a single DM presented in Section \ref{example} and shown in Figure \ref{AttributiTempo}.

%%%%%%%%%%%%%%%%%%%%%%%%%%%%%%%%%%%%%%%
\section{Conclusion}\label{concl}
%%%%%%%%%%%%%%%%%%%%%%%%%%%%%%%%%%%%%%%
In this paper we proposed a general model for combinatorial optimization problems that is based on variables $x_{ilt}$ which take value 1 if facility $i$ is activated in location $l$ at time $t$, and $0$ otherwise. We believe that the model we are proposing has two main merits:
\begin{itemize}
 \item from a more theoretical point of view, our model is in the crossroad of the three following main combinatorial optimization problems:
 \begin{itemize}
    \item  \textit{knapsack problems}, because  our model helps to choose the facilities to be activated as well as the knapsack algorithms determine the items to be selected,
  \item \textit{location problems}, because our model suggests where the selected facilities have to be activated,
  \item \textit{scheduling problems}, because our model suggests also when activating the selected facilities, possibly taking into account some precedence constraints;
 \end{itemize}
  \item from a more application oriented point of view, our model permits to handle complex urban and territorial planning problems in a  multiobjective perspective, taking into account a plurality of stakeholders and policy makers, considering also the uncertainty related to the outcomes of the decision to be taken.
\end{itemize}
Let us point out that our model not necessarily has to be applied to optimization problems with a combinatorial nature. Indeed, for example, the variable $x_{ilt}$ can assume also the meaning of capital allocated to  facility $i$ in location $l$ at time $t$. Therefore the most distinctive feature of our approach is the simultaneous consideration of time and space, so that we refer to our model in terms of space-time model.
With respect to future developments of the research related to the model we are proposing,  the two following points seem to us the most promising:
\begin{itemize}
 \item efficient exact, approximate or heuristic algorithms and procedures to handle problems of big dimensions with many facilities, many constraints and many locations,
 \item applications to real world decision problems in order to test the contribution that our model can give in terms of decision support and to define its possible areas of improvement and enhancement.
\end{itemize}

%%%%%%%%%%%%%%%%%%%%%%%%%%%%
\section*{Acknowledgment}%%%
%%%%%%%%%%%%%%%%%%%%%%%%%%%%
\noindent The research of the second author has been conducted under the Isambard Kingdom Brunel Fellowship Scheme at the University of Portsmouth. The second and the third authors wish to acknowledge funding by the FIR of the University of Catania BCAEA3 “New developments in Multiple Criteria Decision Aiding (MCDA) and their application to territorial competitiveness”. Salvatore Greco has also benefited of the fund “Chance” of the University of Catania.
%\nocite{*}

\textit{}\end{document}